\documentclass[12pt,reqno]{amsart}

\usepackage{amssymb,latexsym}
\usepackage{amsmath}
\usepackage{amsthm}
\usepackage{graphicx}
\usepackage{hyperref}
\usepackage{titletoc}
\usepackage{pdfsync}
\usepackage{color,xcolor}
\usepackage{amsthm,amsmath,amssymb}

\hypersetup{
    colorlinks,%
    citecolor=blue,%
    filecolor=blue,%
    linkcolor=blue,%
    urlcolor=black
}

\usepackage{mathrsfs}

\numberwithin{equation}{section}
\newtheorem{theorem}{Theorem}[section]
\newtheorem{proposition}[theorem]{Proposition}
\newtheorem{lemma}[theorem]{Lemma}
\newtheorem{corollary}[theorem]{Corollary}

\newtheorem*{thmA}{Theorem A}

\usepackage{multirow}
\usepackage[font=bf,aboveskip=15pt]{caption}
\usepackage[toc,page]{appendix}

\theoremstyle{definition}

\newtheorem{remark}[theorem]{Remark}

\def\R{{\mathfrak R}}

\def\begeq{\begin{equation}}
	\def\endeq{\end{equation}}

\def\R{\Bbb R}

\newcommand{\ud}{\mathrm{d}}

\allowdisplaybreaks

\begin{document}

\title[critical Grushin-type problem]{Non-degeneracy and new type of cylindrial solutions for a critical Grushin-type problem}

\author{Yuan Gao, Yuxia Guo and Ning Zhou}

\address{Yuan Gao,
		\newline\indent Department of Mathematical Sciences, Tsinghua University,
		\newline\indent Beijing 100084,  P. R. China.
	}
	\email{gaoy22@mails.tsinghua.edu.cn}

	\address{Yuxia Guo,
		\newline\indent Department of Mathematical Sciences, Tsinghua University,
		\newline\indent Beijing 100084,  P. R. China.
	}
	\email{yguo@tsinghua.edu.cn}

\address{Ning Zhou,
		\newline\indent Department of Mathematical Sciences, Tsinghua University,
		\newline\indent Beijing 100084,  P. R. China.
	}
	\email{zhouning@tsinghua.edu.cn}

\begin{abstract}

In this paper, we consider a critical Grushin-type problem, which is closely related to the prescribed Webster scalar curvature problems on the CR sphere with cylindrically symmetric curvature. We first prove a non-degeneracy result through local Pohozaev identities, then  by using the Lyapunov-Schmidt reduction methods, we construct new type of multi-bubbling solutions with cylindrical symmetry.

\vspace{2mm}
		
		
\vspace{2mm}
		
		
\end{abstract}

\maketitle

\section{Introduction}\label{sec1}

Let  $(\mathbb{S}^{2n+1}, \theta_0)$ be a compact strictly pseudoconvex CR manifold of real dimension $2 n+1$ with the standard contact form $\theta_0$. Given a smooth function $\bar{R}$ on $\mathbb{S}^{2 n+1}$, the prescribed Webster scalar curvature problem on $\mathbb{S}^{2n+1}$ is to find a contact form $\theta$ on $\mathbb{S}^{2n+1}$ conformal equivalent to $\theta_0$ such that the corresponding Webster scalar curvature is $\bar{R}$. If we set $\theta=v^{2 / n} \theta_0$, where $v$ is a smooth positive function on $\mathbb{S}^{2n+1}$, then the above problem is equivalent to solve the following problem:
\begin{equation}\label{eq:S}
    -\Big(2+\frac{2}{n}\Big) \Delta_{\theta_0} v+R_{\theta_0} v=\bar{R} v^{1+\frac{2}{n}}  \quad \text{ on }\, \mathbb{S}^{2n+1},
\end{equation}
where $\Delta_{\theta_0}$ is the sub-Laplacian on $(\mathbb{S}^{2 n+1}, \theta_0)$ and $R_{\theta_0}=n(n+1) / 2$ is the Webster scalar curvature of $(\mathbb{S}^{2 n+1}, \theta_0)$.

Let $\mathbb{H}^n=\mathbb{C}^n \times \mathbb{R} \equiv \mathbb{R}^{2 n+1}$ be the Heisenberg group, using the CR equivalence $F$ (given by the Cayley transform, see \cite{JerisonLeeThe1987}) between $\mathbb{S}^{2 n+1}$ minus a point and $\mathbb{H}^n$, then \eqref{eq:S} becomes (up to an uninfluent constant)
\begin{equation}\label{eq:H}
    -\Delta_{\mathbb{H}^n} u=R u^{\frac{Q+2}{Q-2}} \quad \text{ in } \, \mathbb{H}^n,
\end{equation}
where $\Delta_{\mathbb{H}^n}$ is the canonical sub-elliptic Laplacian on $\mathbb{H}^n$, $Q=2 n+2$ is the homogeneous dimension of $\mathbb{H}^n$, and $R=\bar{R} \circ F^{-1}$. The prescribed Webster scalar curvature problem has been extensively investigated, and many interesting results have been obtained. See, for example, \cite{MalchiodiUguzzoniA2002, FelliUguzzoniSome2004, SalemGamaraThe2011, RiahiGamaraMultiplicity2012, CaoPengYanOn2013, HoPrescribed2016, GamaraHafassaMakniBeta2020, HoKimCR2020} and the references therein.

Denoting by $(Z, t)=(x+i y, t) \equiv(x, y, t)$ the points of $\mathbb{H}^n=\mathbb{C}^n \times \mathbb{R} \equiv \mathbb{R}^{2 n+1}$, we assume that the prescribed curvature $R$ has a natural cylindrical-type symmetry, namely $R(Z, t)=R(|Z|, t)$, which is an analogous case to the radial one in the Euclidean setting.

We will show that under cylindrical-type symmetry assumption, \eqref{eq:H} can be transformed into a Grushin-type equation. The sub-ellipitic Laplacian $\Delta_{\mathbb{H}^n}$ is the second-order differential operator defined as
$$
\Delta_{\mathbb{H}^n}:=\sum_{i=1}^n(X_i^2+Y_i^2),
$$
where
$$
X_i=\frac{\partial}{\partial x_i}+2 y_i \frac{\partial}{\partial t}, \quad Y_i=\frac{\partial}{\partial y_i}-2 x_i \frac{\partial}{\partial t}, \quad i=1, \cdots, n.
$$
Then by direct calculation it holds that
$$
X_i^2 u=\frac{\partial^2 u}{\partial x_i^2}+4 y_i \frac{\partial^2 u}{\partial x_i \partial t}+4 y_i^2 \frac{\partial^2 u}{\partial t^2}, \quad Y_i^2 u=\frac{\partial^2 u}{\partial y_i^2}-4 x_i \frac{\partial^2 u}{\partial y_i \partial t}+4 x_i^2 \frac{\partial^2 u}{\partial t^2} .
$$
Therefore, if $u(Z, t)=u(|Z|, t)>0$ is cylindrical symmetric, problem \eqref{eq:H} becomes
\begin{equation}\label{eq:H-cylindrical}
-\Delta_Z u(|Z|,t)-4|Z|^2 u_{t t}(|Z|,t)=R(|Z|, t) u(|Z|, t)^{\frac{Q+2}{Q-2}}, \quad (Z, t) \in \mathbb{R}^{2 n} \times \mathbb{R},
\end{equation}
where $\Delta_Z$ is the Euclidean Laplacian in $\mathbb{R}^{2 n}$.

\eqref{eq:H-cylindrical} is a special form of the following Grushin-type equation
\begin{equation}\label{eq:Grushin}
    -\Delta_y u(y,z)-4|y|^2 \Delta_z u(y,z)=R(y, z) u(y, z)^{\frac{m_1+2 m_2+2}{m_1+2 m_2-2}},\quad(y, z) \in \mathbb{R}^{m_1} \times \mathbb{R}^{m_2}.
\end{equation}
If $u=u(|y|, z)$ and $R=R(|y|, z)$ satisfy problem \eqref{eq:Grushin}, then we have
$$
-u_{r r}(r, z)-\frac{m_1-1}{r} u_r(r, z)-4 r^2 \Delta_z u(r, z)=R(r, z) u(r, z)^{\frac{m_1+2 m_2+2}{m_1+2 m_2-2}},
$$
where $r=|y|$. Define $v(r, z)=u(\sqrt{r}, z)$, then $v$ satisfies
$$
-v_{r r}(r, z)-\frac{m_1}{2 r} v_r(r, z)-\Delta_z v(r, z)=\frac{R(\sqrt{r}, z)}{4 r} v(r, z)^{\frac{m_1+2 m_2+2}{m_1+2 m_2-2}},
$$
that is, $v=v(|y|, z)>0$ solves the Hardy-Sobolev-type equation
\begin{equation}\label{eq:Hardy-Sobolev}
-\Delta v(y, z)=K(y, z) \frac{v^{\frac{k+h}{k+h-2}}}{|y|}, \quad(y, z) \in \mathbb{R}^k \times \mathbb{R}^h,
\end{equation}
where $k=({m_1+2})/{2}$, $h=m_2$, and $K=K(|y|, z)={R(\sqrt{r}, z)}/{4}$.

A more general Grushin-type equation is
\begin{equation}\label{eq:generalG}
-\Delta_y u(y,z)-(\alpha+1)^2|y|^{2 \alpha} \Delta_z u(y,z)=R(y, z)u(y,z)^{\frac{m_1+(\alpha+1)m_2+2}{m_1+(\alpha+1)m_2-2}}, \quad(y, z) \in\mathbb{R}^{m_1} \times \mathbb{R}^{m_2}.
\end{equation}
Where, the partial differential operator $\mathscr{L}:=\Delta_y+(\alpha+1)^2|y|^{2 \alpha} \Delta_z$ is known as the Grushin operator. The power $\frac{{m_1+(\alpha+1)m_2+2}}{{m_1+(\alpha+1)m_2-2}}$ is the corresponding critical exponent. The Grushin operator is closely related to  the semilinear equations with geometric relevance at the boundary of weakly pseudoconvex domains. Let $\Omega_p=\{(z_1, z_2) \in \mathbb{C}^2 : \operatorname{Im}(z_2)>|z_1|^{2 p}\}$ with $p>1$ be the generalized Siegel domain, which is  a typical example of weakly pseudoconvex domain in the complex space. Under a radial assumption in the variable $z_1$, the natural boundary sub-Laplacian on $\partial \Omega_p$ is the Grushin operator with $\alpha>1$. For more recent results involving the Grushin operator, we refer to \cite{WangWangYangOn2015, LoiudiceAsymptotic2019, LiuTangWangInfinitely2020, LiuNiuConstruction2022, AlvesGandalLoiudiceTyagiA2024} and the references therein.

If $\alpha=0$ and $m_1+m_2=n$, then problem \eqref{eq:generalG} is reduced to
\begin{equation}\label{eq:NirenbergR}
-\Delta u(x)=R(x)u(x)^{\frac{n+2}{n-2}},\quad u>0 \quad \text{ in }\, \R^n.
\end{equation}
Via the sterographic projection, \eqref{eq:NirenbergR} is equivalent to the prescribing scalar curvature problem on the standard $n$-sphere $(\mathbb{S}^n,g_0)$ (i.e., the Nirenberg problem):
\begin{equation}\label{eq:NirenbergS}
-\Delta_{g_0} v+c(n) R_0 v=c(n) R v^{\frac{n+2}{n-2}} \quad \text { on }\, \mathbb{S}^n, \quad \text { for } n \geq 3,
\end{equation}
where $\Delta_{g_0}$ denotes the Laplace-Beltrami operator associated with the metric $g_0$, $c(n)=(n-2) /(4(n-1))$, $R_0=n(n-1)$ is the scalar curvature of $g_0$. There have been many papers on the Nirenberg problem, we refer the readers to \cite{BahriCoronThe1991, ChangGurskyYangThe1993, SchoenZhangPrescribed1996, LiPrescribing1995, LiPrescribing1996} and the references therein. For the generalizations of the Nirenberg problem, please refer to \cite{JinLiXiongOn2014, JinLiXiongOn2015, JinLiXiongThe2017} and references therein.

In this paper, we will consider the following equation
\begin{equation}\label{equation}
-\Delta u(x)=\displaystyle K(x)\frac{u(x)^{2^*-1}}{|y|},\quad u>0 \quad \text{ in } \,  \R^n,
\end{equation}
where $x=(y,z)\in \R^{k}\times\R^{n-k}$, $2\leq k \leq n-1$, $2^{*}:={2(n-1)}/{(n-2)}$.

If $K=K(|y|,z)$ is a cylindrical function, problem \eqref{equation} has been studied extensively. By variational methods, Cao, Peng and Yan \cite{CaoPengYanOn2013} constructed  multi-peak solutions to \eqref{equation} which concentrate exactly at two points between which the distance can be very large. By a Lyapunov-Schmidt reduction argument, Wang, Wang and Yang \cite{WangWangYangOn2015} proved the existence of infinitely many positive solutions with cylindrical symmetry, whose energy can be made arbitrarily large. We refer the readers to \cite{{CatinoLiMonticelliRoncoroniAarXiv}, GheraibiaWangYangExistence2019, {JerisonLeeExtremals1988}, LiuTangWangInfinitely2020, LiuNiuConstruction2022, LiuWangCylindrical2023} for other results of the existence of solutions to \eqref{equation}.

It follows from the classification results of the critical Hardy-Sobolev equation (see \cite{CastorinaFabbriManciniSandeepHardy2009}) that
\begin{equation}\label{bubble}
   U_{\zeta, \mu}(x)= c_{n} \Big(\frac{\mu}{(1+\mu|y|)^2+\mu^2|z-\zeta|^2}\Big)^{\frac{n-2}{2}}, \quad c_{n}=((n-2)(k-1))^{\frac{n-2}{2}}
\end{equation}
is the unique solution to
\begin{equation}\label{limequation}
    -\Delta u(x)=\frac{u(x)^{2^*-1}}{|y|}, \quad u>0, \quad  x=(y,z)\in \R^{k}\times\R^{n-k}.
\end{equation}

Moreover, it follows from  \cite{CastorinaFabbriManciniSandeepHardy2009}, we know that $U_{\zeta, \mu}$ is non-degenerate in
\begin{equation}
    D^{1,2}(\R^n):=\Big\{u:\int_{\R^n}|\nabla u|^2 \,\ud x<+\infty,\, \int_{\R^n}\frac{|u(x)|^{2^*-1}}{|y|} \,\ud x<+\infty\Big\}
\end{equation}
endowed with the inner product $(u,v)=\int_{\R^n}\nabla u\nabla v.$ More precisely, the
kernel of the linear operator associated to \eqref{limequation} is spanned by
\begin{equation}\label{criticalequation}
	\Upsilon_{i}(x) = \frac {\partial U_{0,1}(x)}{\partial z_{i}}, \, i = 1,\cdots,n-k,\quad \Upsilon_{n-k+1}(x) = \frac{n-2}{2}U_{0,1}(x)+x\cdot\nabla U_{0,1}(x).
\end{equation}
Meanwhile, these functions can span the set of the solutions to
\begin{equation}\label{linearized}
	-\Delta u(x)- (2^*-1)\frac{U_{0,1}^{2^*-2}}{|y|} u(x)=0,\quad u \in D^{1,2}(\mathbb{R}^n).
\end{equation}

Define
$$
\begin{aligned}
H_s=\Big\{& u: u \in D^{1,2}(\mathbb{R}^n), \, u \text { is even in } z_2, \\
&u(y, r \cos \vartheta, r \sin \vartheta, z^2)=u\Big(y, r \cos\Big(\vartheta+\frac{2 \pi i}{m}\Big), t \sin \Big(\vartheta+\frac{2 \pi i}{m}\Big), z^2\Big)\Big\}
\end{aligned}
$$
and let $m>0$ be a integer,
\begin{equation}\label{eq:zeta-i}
    \zeta_{i}=\Bigl( \bar r\cos\frac{2(i-1)\pi}m, \bar r\sin\frac{2(i-1)\pi}m, \bar z^2\Bigr),\quad i=1,\cdots,m,
\end{equation}
where $z^2, \bar z^2 \in \R^{n-k-2}.$

Assume that $K(x)$ satisfy the following conditions:\\
${\bf (K_1')}$: $K(x)=K(|z^1|,z^2)\ge0$ are bounded functions for $x=(y,z^1,z^2)\in \R^{k}\times\R^{2}\times\R^{n-k-2}$. Set $r:=|z^1|$, $K(r,z^2)$ has a stable critical point $(r_0,z_0^2)$ satisfying $r_0>0$, $K(r_0,z_0^2)=1>0$ and
$$
\deg (\nabla K(r,z^2), (r_0,z_0^2))\neq 0;
$$\\
${\bf (K_2')}$: $K(r,z^2)\in C^3(B_{\rho'_0}(r_0,z_0^2))$ for $\rho'_0>0$ is a fixed small constant, and
$$
\Delta K(r_0,z_0^2):=\frac{\partial^2K}{\partial r^2}(r_0,z_0^2)+\sum_{i=3}^{n-k}\frac{\partial^2K}{\partial z_i^2}(r_0,z_0^2)<0.
$$
Then under the assumptions of ${\bf (K_1^\prime)}$-${\bf (K_2^\prime)}$, Liu and Wang \cite{LiuWangCylindrical2023} obtained the existence of bubble solutions to \eqref{equation}. Their result states as the following:

\begin{thmA}\label{thmA}
Suppose that $n\geq 5$, $\frac{n+1}{2}\le k<n-1$, $K(x)$ satisfies ${\bf (K_1^\prime)}$-${\bf (K_2^\prime)}$, then there exists an integer $m_0>0$, such that for any integer $m>m_0$, problem \eqref{equation} has a solution $u_m$ of the form
\begin{equation}\label{form}
    u_m=\sum_{i=1}^m {\tilde\eta} U_{\zeta_i, \mu_m}+\phi_m,
\end{equation}
where ${\tilde\eta}\in [0,1]$ is some cut-off function such that $\tilde\eta(x)=1$ if $|(|y|,r,z^2)-(0,r_0,z_0^2)|\le\delta'$, and $\tilde\eta(x)=0$ if $|(|y|,r,z^2)-(0,r_0,z_0^2)|\ge 2\delta'$ with $\delta'>0$ a small constant satisfying $K(r,z^2)\ge C>0$ for $|(r,z^2)-(r_0,z_0^2)|\le 10 \delta'$, and $\zeta_i$ is defined as in \eqref{eq:zeta-i}, $\phi_m \in  {H_s}$.

Furthermore, as $m \to+\infty$, $(\bar r_m, \bar z^2_m) \to (r_0, z^2_0)$, $\mu_m \in [L_0' m^{\frac{n-2}{n-4} }$, $L_1' m^{\frac{n-2}{n-4}}]$, $L_1'>L_0'>0$ are some constants, and $\|\phi_m\|_{L^{\infty}(\R^n)}=o(\mu_m^{\frac{n-2}{2}})$.
\end{thmA}

In order to obtain the solutions of the form \eqref{form}, the authors gluing a very large number of basic profiles \eqref{bubble} together, which centered at the vertices of a regular polygon with a large number of edges. Note that the solution $u_m$ is radially in  $z^1$. And  of course,  by the same argument, we can also construct a solution $u_q$ with $q$-bubbles, and $u_q$ is radially with respect to the first two components of $z^2$.

In this paper, we want to discuss whether $u_m$ and $u_q$ can be glued together to give rise to a new type of solutions, with $m$ and $q$ possibly being different orders. Specifically, we want to construct a new solution to \eqref{equation} whose shape is, at main order,
\begin{equation}\label{eq:new}
    u \approx \sum_{i=1}^m {\eta}U_{\tilde\zeta_i, \tilde\mu_m}+\sum_{j=1}^q \eta U_{p_j, \lambda_q},
\end{equation}
for $m$ and $q$ large,
where $\eta$ are some cut-off functions defined later, and
\begin{equation}\label{eqtildezeta-i}
\tilde\zeta_i=\Bigl( \bar r\cos\frac{2(i-1)\pi}m, \bar r\sin\frac{2(i-1)\pi}m, 0, 0, \tilde z'\Bigr),\quad \tilde z'\in \R^{n-k-4},\quad i=1,\cdots,m,
\end{equation}
\begin{equation}\label{eq:p-j}
p_j=\Bigl( 0,0,\bar t\cos\frac{2(j-1)\pi}q, \bar t\sin\frac{2(j-1)\pi}q, \bar z'\Bigr),\quad\bar z'\in \R^{n-k-4},\quad j=1,\cdots,q.
\end{equation}

Notice that it's very difficult to obtain solution \eqref{eq:new} by perturbation arguments. In fact, if we want to make a small correction to obtain a solution to \eqref{equation} of shape \eqref{eq:new} with $q \gg m$, the estimate of the correction term is dominated by the parameter $m$. In other words, it is hard to see the contribution to the energy from the bumps $U_{p_j,\lambda_q}$. Therefore, it is not easy to directly to construct solutions of the form \eqref{eq:new}.

To overcome this difficulty, we use a new method which was first introduced by Guo, Musso, Peng and Yan in \cite{GuoMussoPengYanNon2020}. They first proved a non-degeneracy result for the positive multi-bubbling solutions of the prescribed scalar curvature equations constructed in \cite{WeiYanInfinitely2010}. Then they used this non-degeneracy result to glue together bubbles with different concentration rate to obtain new solutions. We refer to \cite{GuoMussoPengYanNon2022, NieLiNon2022, HeWangWangNew2022, GuoHuLiuNon2023, GaoGuoHuNon2024} for the applications of this method to various problems.

In order to obtain the non-degeneracy result, we assume further \\
${\bf (K_3')}$: The matrix
\begin{equation}
    \left(\begin{matrix}
    \displaystyle\frac{\partial^2 K(r_0,z_0^2)}{\partial z_1^2}&\quad \displaystyle\frac{\partial^2 K(r_0,z_0^2)}{\partial z_1\partial z_3} \quad \cdots \quad \frac{\partial^2 K(r_0,z_0^2)}{\partial z_1\partial z_{n-k}}\\
    \\
    \displaystyle\frac{\partial^2 K(r_0,z_0^2)}{\partial z_3\partial z_1}&\quad \displaystyle\frac{\partial^2 K(r_0,z_0^2)}{\partial z_3^2} \quad \cdots \quad \frac{\partial^2 K(r_0,z_0^2)}{\partial z_3\partial z_{n-k}}\\
    \\
    \vdots&\\
    \\
    \displaystyle\frac{\partial^2 K(r_0,z_0^2)}{\partial z_{n-k}\partial z_1}&\quad \displaystyle\frac{\partial^2 K(r_0,z_0^2)}{\partial z_{n-k}\partial z_3} \quad \cdots \quad \frac{\partial^2 K(r_0,z_0^2)}{\partial z_{n-k}^2}\\
\end{matrix}\right)
\end{equation}
is non-degenerate.\\

Here is our first result:
\begin{theorem}\label{nondeg}
    Suppose that $n\ge 8$, $\frac{n+1}{2}\le k<n-3$, $K(x)$ satisfies ${\bf (K_1')}$-${\bf (K_3')}$, then there exists a large $m_0,$ such that for any integer $m>m_0,$  the bubble solution in Theorem A is non-degenerate, in the sense that if $\xi\in {H_s}$ is a solution of the following linear equation:
$$
    L_{m}\xi:=-\Delta\xi-(2^*-1)K(x)\frac{u_m^{2^*-2}}{|y|}\xi=0 \quad \text{ in } \, \R^n,
$$
then $\xi=0$.
\end{theorem}

As an application of Theorem \ref{nondeg}, we define the symmetric Sobolev space
\begin{equation}\label{eq:Xs}
\begin{aligned}
X_s=\Big\{&u: u \in D^{1,2}(\mathbb{R}^n),\, u \text { is even in } z_h,\, h=1,2,4,  \\
&u(y, z_1, z_2, t \cos \theta, t \sin \theta, z^{\prime})=u\Big(y, z_1, z_2, t \cos \Big(\theta+\frac{2 \pi j}{q}\Big), t \sin \Big(\theta+\frac{2 \pi j}{q}\Big), z^{\prime}\Big)\Big\}.
\end{aligned}
\end{equation}
Since we aim to glue the bubbles centered at $(z_1, z_2)$-plane and $(z_3, z_4)$-plane separately, the main term $\sum_{i=1}^m{\eta}U_{\tilde\zeta_i, \tilde\mu_m}+\sum_{j=1}^q \eta U_{p_j, \lambda_q}$  is in the symmetric Sobolev space $H_s\cap X_s$. Note that we will construct solution which may radially in  depend on $z^*=(z_1,z_2,z_3,z_4)$. So we improve our assumptions on $K(x)$ correspondingly. More precisely, we assume $K(x)$ satisfies:\\
${\bf (K_1)}$: $K(x)=K(|z^*|,z')\ge0$ is a bounded function for $x=(y,z^*,z')\in \R^{k}\times\R^{4}\times\R^{n-k-4}$. Set $t:=|z^*|$, $K(t,z')$ has a stable critical point $(t_0,z_0')$ satisfying $t_0>0$, $K(t_0,z_0')=1$ and
$$
\deg (\nabla K(t,z'), (t_0,z_0'))\neq 0;
$$
${\bf (K_2)}$: $K(t,z')\in C^3(B_{\rho_0}(t_0,z_0'))$, where $\rho_0>0$ is a fixed small constant, and
$$
\Delta K(t_0,z_0')<0;
$$
${\bf (K_3)}$: The matrix
\begin{equation}
    \left(\begin{matrix}
    \displaystyle\frac{\partial^2 K(t_0,z_0')}{\partial z_1^2}&\quad \displaystyle\frac{\partial^2 K(t_0,z_0')}{\partial z_1\partial z_3} \quad \cdots \quad \frac{\partial^2 K(t_0,z_0')}{\partial z_1\partial z_{n-k}}\\
    \\
    \displaystyle\frac{\partial^2 K(t_0,z_0')}{\partial z_3\partial z_1}&\quad \displaystyle\frac{\partial^2 K(t_0,z_0')}{\partial z_3^2} \quad \cdots \quad \frac{\partial^2 K(t_0,z_0')}{\partial z_3\partial z_{n-k}}\\
    \\
    \vdots&\\
    \\
    \displaystyle\frac{\partial^2 K(t_0,z_0')}{\partial z_{n-k}\partial z_1}&\quad \displaystyle\frac{\partial^2 K(t_0,z_0')}{\partial z_{n-k}\partial z_3} \quad \cdots \quad \frac{\partial^2 K(t_0,z_0')}{\partial z_{n-k}^2}\\
\end{matrix}\right)
\end{equation}
is non-degenerate.

\begin{remark}\label{remark1}
According to the proof of Theorem A in \cite{LiuWangCylindrical2023}, under the condition ${\bf (K_1)}$-${\bf (K_2)}$ instead of ${\bf (K_1^{\prime})}$-${\bf (K_2^{\prime})}$, we can also proof the result of Theorem A is true by a similar argument. That is, we can obtain bubble solutions to \eqref{equation} centered at the point $\tilde\zeta_i$. For simplicity of notation, we still denote the solution as $u_m$.  We leave some detail of $u_m$ in Section \ref{constru}. Moreover, according to proof of Theorem \ref{nondeg}, we can similarly deduce that under ${\bf (K_1)}$-${\bf (K_3)}$, the solution $u_m$ mentioned above is also non-degenerate. We denote $\tilde\zeta_i$, $\tilde \mu_m$ as $\zeta_i$, $\mu_m$ without causing ambiguity when constructing new kinds of solutions in  Theorem \ref{newsolution}.
\end{remark}

Let $\bar\delta>0$ be a small constant satisfying $K(t,z')\ge C>0$ for $|(t,z')-(t_0,z_0')|\le 10\bar\delta$. We define a cut-off function $\eta(x)=\eta(|y|,|z^*|,z')\in[0,1]$ such that $\eta(x)=1$ if $|(|y|,t,z')-(0,t_0,z_0')|\le\bar\delta$, and $\eta(x)=0$ if $|(|y|,t,z')-(0,t_0,z_0')|\ge 2\bar\delta$.
{{We always assume that $|(\bar t, \bar z^2)-(r_0, z^2_0)|\le\frac{1}{\mu^{1-\bar\theta}}$ with some constant $\bar \theta\in(0,{(1-\epsilon_0)}/{2})$ and ${n}/{2}-\bar\theta-\tau>2$, for $\epsilon_0>0$ is a fixed constant taken later in Lemma \ref{lemma2.3}.}}

We have
\begin{theorem}\label{newsolution}
   Suppose that $n\ge 8$, $\frac{n+1}{2}\le k<n-3$, $K(x)$ satisfies ${\bf (K_1)}$-${\bf (K_3)}$,  and assume that $u_m$ is the solution to \eqref{equation} gotten from Remark \ref{remark1} with $m>0$ a large even integer, then there exists an integer $q_0>0$, such that for any integer $q>q_0$, problem \eqref{equation} has a solution $v_q$ of the form
\begin{equation}\label{eq:newsolution}
v_q=u_m+\sum_{j=1}^{q} \eta U_{p_j,\lambda_q}+\psi_q,
\end{equation}
where $p_j$ is defined as in \eqref{eq:p-j}, $\psi_q \in  {X_s}$, $(\bar t_q, \bar z'_q) \to (t_0, z'_0)$, $\lambda_q \in [L_0 q^{\frac{n-2}{n-4} }$, $L_1 q^{\frac{n-2}{n-4}}]$, $L_1>L_0>0$ are some constants, and $\|\psi_q\|_{L^{\infty}(\R^n)}=o(\lambda_q^{\frac{n-2}{2}}).$
\end{theorem}

As a result of Theorem \ref{newsolution} and the equivalence of equations \eqref{eq:Grushin} and \eqref{eq:Hardy-Sobolev}, we can obtain the existence of cylindrically symmetric multi-bubbling solutions to the critical Grushin-type equation \eqref{eq:Grushin}:

\begin{corollary}\label{cor:Grushin}
Assume that $R(y, z)=R(|y|, z)$ is bounded and continuous in $\mathbb{R}^{m_1+m_2}$. Also assume that $K=K(|y|, z)={R(\sqrt{|y|}, z)}/{4}$ satisfies ${\bf (K_1)}$-${\bf (K_3)}$. Then problem \eqref{eq:Grushin} has infinitely many cylindrically symmetric multi-bubbling solutions.
\end{corollary}

This paper is organized as follows. In Section 2 we prove the non-degeneracy result stated in Theorem \ref{nondeg}, which is a important ingredient in the construction of the new type of bubbling solutions. Using this non-degeneracy result, we prove Theorem \ref{newsolution} in Section 3. We present some important identities and essential estimates,  which are used in Sections 2 and 3, in Appendices.

\section{Non-degeneracy of the bubbling solutions}
In this section, we will prove the non-degeneracy of the multi-bubbling solutions obtained in Theorem A. Let us first introduce the following weighted norms:
$$
    \|u\|_{ *}:= \sup_{x\in \R^n}\Big(\sum_{j=1}^{m}\frac{\mu^{\frac{n-2}{2}}}{(1+\mu  |y|+\mu |z-\zeta_j|)^{\frac{n-2}{2}+\tau}}\Big)^{-1}|u(x)|,
$$
$$
    \|f\|_{ {**}}:=\sup_{x\in \R^n}\Big(\sum_{j=1}^{m}\frac{\mu^{\frac{n+2}{2}}}{\mu |y|(1+\mu  |y|+\mu |z-\zeta_j|)^{\frac{n}{2}+\tau}}\Big)^{-1}|f(x)|,
$$
where $\tau=\frac{n-4}{n-2}$.
 Denote
$$
\overline{U}_{\zeta_j,\mu}=\tilde\eta(x) U_{\zeta_j,\mu}, \  \overline{W}_{\bar r,\bar z^2,\mu}=\sum_{j=1}^{m}\overline U_{\zeta_j,\mu}, \ W_{\bar r,\bar z^2,\mu}=\sum_{j=1}^{m}U_{\zeta_j,\mu},
$$
where $\tilde\eta$ is as in Theorem A.
Throughout our paper, we employ  $\delta$, $\epsilon$, $\varepsilon$, $\sigma$ to denote some small constants.

\begin{lemma}\label{lemma2.1}
There exists a constant $C>0$ such that
\begin{equation}\label{2-1}
|u_m(x)|\le C\sum_{j=1}^{m} \frac{\mu^{\frac{n-2}{2}}}{(1+\mu|y|+\mu|z-\zeta_j|)^{n-2}}.
\end{equation}
\end{lemma}

\begin{proof}
By Green's representation, H\"older inequality and Lemma \ref{B.6}, we have that
\begin{align*}
    &|u_m(x)|\leq C\int_{\R^n}\frac{1}{|\tilde x-x|^{n-2}}K(\tilde r,\tilde z^2)\frac{u_m^{2^*-1}}{|\tilde y|} \,\text{d}\tilde x\\
    \le&C\int_{\R^n} \frac{1}{|\tilde x-x|^{n-2}|\tilde y|}\Big(\sum_{j=1}^{m}\frac{\mu^{\frac{n-2}{2}}}{(1+\mu|\tilde y|+\mu|\tilde z-\zeta_j|)^{n-2}}\Big)^{2^*-1}\\
    &+C\|\phi_m\|_{*}^{2^*-1}\int_{\R^n}  \frac{1}{|\tilde x-x|^{n-2}|\tilde y|}\Big(\sum_{j=1}^{m}\frac{\mu^{\frac{n-2}{2}}}{(1+\mu|\tilde y|+\mu|\tilde z-\zeta_j|)^{\frac{n-2}{2}+\tau}}\Big)^{2^*-1}\\
    \le&C\int_{\R^n} \frac{1}{|\tilde x-x|^{n-2}}\frac{1}{\mu|\tilde y|}\sum_{j=1}^{m}\frac{\mu^{\frac{n+2}{2}}}{(1+\mu|\tilde y|+\mu|\tilde z-\zeta_j|)^{n-\frac{2}{n-2}\tau}}\Big(1+\sum_{j=2}^{m}\frac{1}{(\mu|\zeta_j-\zeta_1|)^\tau}\Big)^{\frac{2}{n-2}}\\
    &+C\frac{1}{\mu^{(1+\epsilon)(2^*-1)}}\int_{\R^n} \frac{1}{|\tilde x-x|^{n-2}}\frac{1}{\mu|\tilde y|}\sum_{j=1}^{m}\frac{\mu^{\frac{n+2}{2}}}{(1+\mu|\tilde y|+\mu|\tilde z-\zeta_j|)^{\frac{n}{2}+\frac{n}{n-2}\tau-\frac{2}{n-2}\tau_1}}\\
    &\times\Big(1+\sum_{j=2}^{m}\frac{1}{(\mu|\zeta_j-\zeta_1|)^{\tau_1}}\Big)^{\frac{2}{n-2}}\\
    \le&C\sum_{j=1}^{m}\frac{\mu^{\frac{n-2}{2}}}{(1+\mu|y|+\mu|z-\zeta_j|)^{n-1-\frac{2}{n-2}\tau}}+C\sum_{j=1}^{m}\frac{\mu^{\frac{n-2}{2}}}{(1+\mu|y|+\mu|z-\zeta_j|)^{\frac{n-2}{2}+\frac{n}{n-2}\tau-\frac{2}{n-2}\tau_1}}
\end{align*}
for $0<\tau_1<\tau$. Since that
$$
n-1-\frac{2}{n-2}\tau>n-2,
$$
and
$$
\frac{n-2}{2}+\frac{n}{n-2}\tau-\frac{2}{n-2}\tau_1=\frac{n-2}{2}+\tau+\frac{2}{n-2}(\tau-\tau_1)>\frac{n-2}{2}+\tau,
$$
then we can continue this process and finally obtain \eqref{2-1}.
\end{proof}

In the following, we will apply local Pohozaev identities to prove the non-degeneracy of the bubbling solutions. We argue by contradiction.
Suppose that there exists $m_\ell\to +\infty,$ satisfying
$$
L_{m_\ell}\xi_\ell=0 \quad  \text{ in }\, \mathbb{R}^{n},
$$
but $\xi_\ell\not\equiv0.$ Without loss of generality, we may assume $\|\xi_{\ell}\|_{*} = 1$ and obtain the contradictions through the following steps.  Define
$$
\widehat{\xi}_{\ell}(x) = \mu_{m_\ell}^{-\frac{n-2}{2}}\xi_{\ell}(\mu_{m_\ell}^{-1}x+(0,\zeta_1)),
$$
where $\zeta_1$ is as in \eqref{eq:zeta-i}.
\begin{lemma}\label{lemma2.2}
It holds
$$
\widehat{\xi}_{\ell} \to b_{0}\Phi_{0} + b_{1}\Phi_{1}+\sum_{i=3}^{n-k}b_{i}\Phi_{i},
$$
uniformly in $C^{1}(B_{R}(0))$ for any $R > 0$, where $b_{0}$ and $b_{i}$, $i=1, 3, 4,\cdots n-k$ are some constants,
$$
\Phi_{0} = \frac{\partial U_{0,\mu}}{\partial \mu}\Big|_{\mu = 1},\quad  \Phi_{i} = \frac{\partial U_{0,1}}{\partial z_{i}}, \, i = 1,\cdots,n-k.
$$
\end{lemma}

\begin{proof}
By $\|\xi_\ell\|_{*}=1$, we have $|\widehat{\xi}_{\ell}| \leq C$. Therefore, we may assume that $\widehat{\xi}_{\ell} \to \xi$ in $C^1_{loc}(\mathbb{R}^{n})$. Then $\xi$ satisfies
$$
-\Delta \xi = (2^{*}-1)\frac{U_{0,1}^{2^*-2}}{|y|}\xi \quad \text{ in } \,  \mathbb{R}^{n},
$$
which gives
$$
\xi = \sum_{i=0}^{n-k}b_{i}\Phi_{i}.
$$
Since $\xi_{\ell}$ is even in $z_{2}$, it holds that $b_{2} = 0$.
\end{proof}

We decompose
$$
\xi_{\ell}(x) = b_{0,\ell}\mu_{{m_\ell}}\sum_{j=1}^{{m_\ell}}\frac{\partial \overline{U} _{\zeta_j,\mu_{{m_\ell}}}}{\partial \mu_{{m_\ell}}}
    - b_{1,\ell}\mu_{{m_\ell}}^{-1}\sum_{j=1}^{{m_\ell}}\frac{\partial  \overline{U}_{\zeta_j,\mu_{{m_\ell}}}}{\partial \bar r}
	-\sum_{i=3}^{n-k} b_{i,\ell}\mu_{{m_\ell}}^{-1}\sum_{j=1}^{{m_\ell}}\frac{\partial  \overline{U}_{\zeta_j,\mu_{{m_\ell}}}}{\partial \bar z_i}  + \xi_{\ell}^{*},
$$
where
$\xi_{\ell}^{*}$ satisfies that, for $i=3,\cdots,n-k,$
$$
\int_{\mathbb{R}^{n}}\frac{\overline{U} _{\zeta_j,\mu_{{m_\ell}}}^{2^{*}-2}}{|y|} \frac{\partial \overline{U} _{\zeta_j,\mu_{{m_\ell}}}}{\partial \mu_{{m_\ell}}}\xi_{\ell}^{*} = \int_{\mathbb{R}^{n}} \frac{\overline{U} _{\zeta_j,\mu_{{m_\ell}}}^{2^{*}-2}}{|y|} \frac{\partial \overline{U} _{\zeta_j,\mu_{{m_\ell}}}}{\partial \bar r}\xi_{\ell}^{*}
		=\int_{\mathbb{R}^{n}} \frac{\overline{U} _{\zeta_j,\mu_{{m_\ell}}}^{2^{*}-2}}{|y|} \frac{\partial \overline{U} _{\zeta_j,\mu_{{m_\ell}}}}{\partial \bar z_i}\xi_{\ell}^{*}  = 0.
$$
It follows from Lemma \ref{lemma2.2} that $b_{i,\ell}$ are bounded for $ i=1,3,\cdots,n-k$. We first give an estimate to $\xi_{\ell}^{*} $.

\begin{lemma}\label{lemma2.3}
It holds
\begin{align}\label{xi_n}
\|\xi_{\ell}^{*}\|_{*} \leq  \frac{C}{\mu_{m_\ell}^{\frac{n-2\tau}{n-2}} }.
\end{align}
\end{lemma}

\begin{proof}
A direct calculation leads to that
$$
\begin{aligned}
        L_{m_\ell}\xi_{\ell}^{*}=&-\Delta\xi_{\ell}^{*}-(2^*-1)K(r,z^2)\frac{u^{2^*-2}_{m_\ell}}{|y|}\xi_{\ell}^{*}\\
        =&(2^*-1)\tilde\eta(x)(K(r,z^2)-1)\frac{u^{2^*-2}_{m_\ell}}{|y|}\sum_{j=1}^{m_\ell}  \beta_j+(2^*-1)\tilde\eta(x)\sum_{j=1}^{m_\ell}\Big(\frac{u^{2^*-2}_{m_\ell}}{|y|}-\frac{\overline{U}^{2^*-2} _{\zeta_j,\mu_{{m_\ell}}}}{|y|}\Big)  \beta_j\\
        &+\Delta\tilde\eta(x)\sum_{j=1}^{m_\ell}  \beta_j
        +2\nabla\tilde\eta(x)\sum_{j=1}^{m_\ell}\nabla\beta_j\\
        :=&I_1+I_2+I_3+I_4,
\end{aligned}
$$
where
$$
    \beta_j:=b_{0,\ell}\mu_{{m_\ell}}\frac{\partial \overline{U} _{\zeta_j,\mu_{{m_\ell}}}}{\partial \mu_{{m_\ell}}}
		- b_{1,\ell}\mu_{{m_\ell}}^{-1}\frac{\partial  \overline{U}_{\zeta_j,\mu_{{m_\ell}}}}{\partial \bar r}
		-\sum_{i=3}^{n-k} b_{i,\ell}\mu_{{m_\ell}}^{-1}\frac{\partial  \overline{U}_{\zeta_j,\mu_{{m_\ell}}}}{\partial \bar z_i}.
$$
In the following, we estimate the terms above one by one. Define
$$
\begin{aligned}
    \Omega_j:=\Big\{&x:x=(y,z_1,z_2,z'')\in\R^k\times\R\times\R\times\R^{n-k-2},\\
&\Big\langle\frac{(z_1,z_2)}{|(z_1,z_2)|},\, \Big(\cos\frac{2(j-1)\pi}{m},\sin\frac{2(j-1)\pi}{m}\Big)\Big\rangle\ge\cos\frac{\pi}{m}\Big\}.
\end{aligned}
$$
Without loss of generality, we may assume that $y\in \Omega_1.$ For $I_1,$ we have
\begin{align}\label{I1}
    |I_1|\le& C\frac{|K(r,z^2)-1|}{|y|}\Big(\sum_{j=1}^{m_\ell}\frac{\mu_{{m_\ell}}^{\frac{n-2}{2}}}{(1+\mu_{{m_\ell}}|y|+\mu_{{m_\ell}}|z-\zeta_j|)^{n-2}}\Big)^{2^*-1}\nonumber\\
    \le&
    C\frac{|K(r,z^2)-1|\mu_{{m_\ell}}^{\frac{n+2}{2}}}{\mu_{{m_\ell}}|y|(1+\mu_{{m_\ell}}|y|+\mu_{{m_\ell}}|z-\zeta_1|)^{n}}\nonumber\\
    &+C\frac{\mu_{{m_\ell}}^{\frac{n-2}{2}}}{|y|(1+\mu_{{m_\ell}}|y|+\mu_{{m_\ell}}|z-\zeta_1|)^{n-2}}\Big(\sum_{j=2}^{m_\ell}\frac{\mu_{{m_\ell}}^{\frac{n-2}{2}}}{(1+\mu_{{m_\ell}}|y|+\mu_{{m_\ell}}|z-\zeta_j|)^{n-2}}\Big)^{2^*-2}\nonumber\\
    :=& I_{11}+I_{12}.
\end{align}
By the Taylor expansion to $K(r,z^2)$, for $|(r,z^2)-(r_0,z_0^2)|\le\frac{\delta'}{\mu_{{m_\ell}}^{({1+\epsilon_0})/{2}}}<\rho_0$, where $\epsilon_0>0$ is a small constant fixed later, we have
\begin{equation}\label{I111}
    |I_{11}|\le \frac{C}{\mu_{{m_\ell}}^{1+\epsilon_0}} \frac{\mu_{{m_\ell}}^{\frac{n+2}{2}}}{\mu_{{m_\ell}}|y|(1+\mu_{{m_\ell}}|y|+\mu_{{m_\ell}}|z-\zeta_1|)^{\frac{n}{2}+\tau}}.
\end{equation}
On the other hand, for $\frac{\delta'}{\mu_{{m_\ell}}^{{(1+\epsilon_0)}/{2}}}\le |(r,z^2)-(r_0,z_0^2)|\le \delta',$ we have $$|(r,z^2)-(\bar r,\bar z^2)|\ge \frac{\delta'}{\mu_{{m_\ell}}^{{(1+\epsilon_0)}/{2}}}-\frac{1}{\mu_{{m_\ell}}^{1-\bar\theta}}\ge\frac{\delta'}{2\mu_{{m_\ell}}^{({1+\epsilon_0})/{2}}},$$ since $\bar\theta<\frac{1-\epsilon_0}{2}.$ Then
\begin{equation}\label{I112}
|I_{11}|\le C\frac{1}{\mu_{{m_\ell}}^{\frac{1-\epsilon_0}{2}(\frac{n}{2}-\tau)}}\sum_{j=1}^{m}\frac{\mu_{{m_\ell}}^{\frac{n+2}{2}}}{\mu_{{m_\ell}}|y|(1 +\mu_{{m_\ell}}|y|+ \mu_{{m_\ell}}|z-{\zeta_j}|)^{\frac{n}{2}+\tau}}.
\end{equation}
By \eqref{I111} and \eqref{I112}, we have
\begin{align}\label{I11}
\|I_{11}\|_{**}\le \frac{1}{\mu_{{m_\ell}}^{\min\{1+\epsilon_0,\frac{1-\epsilon_0}{2}(\frac{n}{2}-\tau)\}}}.
\end{align}
For $I_{12},$ we can check that
\begin{equation}\label{I12}
\begin{aligned}
        |I_{12}|\le&\Big(\sum_{j=2}^{m_\ell}\frac{1}{\mu_{{m_\ell}}^{\frac{n}{2}-\tau}|{\zeta_j}-{\zeta_1}|^{\frac{n}{2}-\tau}}\Big)\frac{\mu_{{m_\ell}}^{\frac{n+2}{2}}}{\mu_{{m_\ell}}|y|(1 +\mu_{{m_\ell}}|y|+ \mu_{{m_\ell}}|z-{\zeta_1}|)^{\frac{n}{2}+\tau}}\\
        \le&C\Big(\frac{m_\ell}{\mu_{{m_\ell}}}\Big)^{\frac{n}{2}-\tau}\frac{\mu_{{m_\ell}}^{\frac{n+2}{2}}}{\mu_{{m_\ell}}|y|(1 +\mu_{{m_\ell}}|y|+ \mu_{{m_\ell}}|z-{\zeta_1}|)^{\frac{n}{2}+\tau}}.
\end{aligned}
\end{equation}
Since we can always take a proper $\epsilon_0$ to make
$$
\frac{1}{\mu_{{m_\ell}}^{\min\{1+\epsilon_0, \frac{1-\epsilon_0}{2}(\frac{n}{2}-\tau) \}}}=o\Big(\frac{m_\ell}{\mu_{{m_\ell}}}\Big)^{\frac{n}{2}-\tau},
$$
therefore, combining \eqref{I1}--\eqref{I12}, we finally get
\begin{equation}\label{I_1}
\|I_1\|_{**} \le C\Big(\frac{m_\ell}{\mu_{{m_\ell}}}\Big)^{\frac{n}{2}-\tau}.
\end{equation}

Next, we estimate $I_2$, similar to $I_{12}$, we can easily get
\begin{align}\label{I2}
    |I_2|\le &C\sum_{j=1}^{m_\ell}\Big(\frac{u^{2^*-2}_{m_\ell}}{|y|}-\frac{\overline{U}^{2^*-2} _{\zeta_j,\mu_{{m_\ell}}}}{|y|}\Big){U} _{\zeta_j,\mu_{{m_\ell}}}\nonumber\\
    \le & C\frac{1}{|y|}{U}^{2^*-2} _{\zeta_j,\mu_{{m_\ell}}}\Big(\sum_{j=2}^{m_\ell}{U}^{2^*-2} _{\zeta_j,\mu_{{m_\ell}}}+|\phi_{m_\ell}| \Big)\nonumber\\
     \le&C\Big(\frac{m_\ell}{\mu_{{m_\ell}}}\Big)^{\frac{n}{2}-\tau}\frac{\mu_{{m_\ell}}^{\frac{n+2}{2}}}{\mu_{{m_\ell}}|y|(1 +\mu_{{m_\ell}}|y|+ \mu_{{m_\ell}}|z-{\zeta_1}|)^{\frac{n}{2}+\tau}}.
\end{align}
Then
\begin{equation}\label{I_2}
\|I_2\|_{**} \le C\Big(\frac{m_\ell}{\mu_{{m_\ell}}}\Big)^{\frac{n}{2}-\tau}.
\end{equation}

Noting that for $x\in \operatorname{supp}|\nabla\tilde\eta|$, $1+\mu_{m_\ell}|y|+\mu_{m_\ell}|z-\zeta_{i}|\ge C\mu$, we can get estimates for $I_3$:
\begin{equation}\label{I3}
\begin{aligned}
    |I_3|\le &C\sum_{j=1}^{m_\ell}\frac{\tilde\eta\mu_{m_\ell}^{\frac{n-2}{2}}}{(1 +\mu_{m_\ell}|y|+ \mu_{m_\ell}|z-{\zeta_j}| )^{n-2}}\\
    \le& \frac{C}{\mu_{m_\ell}^{\frac{n-2}{2}-\tau}}\sum_{j=1}^{m_\ell}\frac{\mu_{m_\ell}^{\frac{n+2}{2}}}{\mu_{m_\ell}|y|(1 +\mu_{m_\ell}|y|+ \mu_{m_\ell}|z-{\zeta_j}| )^{\frac{n}{2}+\tau}}.
\end{aligned}
\end{equation}
Thus,
\begin{equation}\label{I_3}
\|I_3\|_{**} \le  \frac{C}{\mu_{m_\ell}^{\frac{n-2}{2}-\tau}},
\end{equation}
and similar to the estimation of $\|I_3\|_{**}$, it also holds that
\begin{equation}\label{I_4}
\|I_4\|_{**} \le  \frac{C}{\mu_{m_\ell}^{\frac{n-2}{2}-\tau}}.
\end{equation}
Combining \eqref{I_1}, \eqref{I_2}, \eqref{I_3}, \eqref{I_4}, and similar to the proof in \cite{LiuWangCylindrical2023}, we can prove that there exist a constant $\varrho>0$ such that
$$
\|\xi_\ell^*\|_{*}\le \frac{1}{\varrho} \|L_{m_\ell}\xi_\ell^*\|_{**}\le \frac{C}{\mu_{m_\ell}^{\min\{\frac{n}{n-2}-\frac{2}{n-2}\tau,\frac{n-2}{2}-\tau \}}}=\frac{C}{\mu_{m_\ell}^{\frac{n}{n-2}-\frac{2}{n-2}\tau }}.
$$
\end{proof}

\begin{proposition}\label{pro2.4}
If ${\bf (K_3)}$  holds, then
    $\widehat{\xi}_{\ell}\to 0$
uniformly in $C^1(B_R(0))$ for any $R>0$.
\end{proposition}

\begin{proof}
The proof consists of the following steps.

{\bf Step 1.} We first prove $b_{i,\ell}\to 0$, $i=1,3,4,\cdots,n-k,$ by applying local Pohozaev identity \eqref{pohozeav1} and \eqref{pohozeav2} in $\Omega_1$. By the symmetry, we have $\frac{\partial u_{m_\ell}}{\partial \nu}=\frac{\partial \xi_\ell}{\partial \nu}=0$ and $\langle\nu, y\rangle=0$ on $\partial\Omega_1$. Then we have
\begin{equation}\label{po1}
- \int_{\Omega_1}\frac{\partial K(r,z^2)}{\partial z_j}\frac{u_{m_\ell}^{2^*-1}\xi_{\ell}}{|y|}=\int_{\partial\Omega_1}\nabla u_{m_\ell}\nabla \xi_{\ell} \nu_{k+j}-\int_{\partial\Omega_1}K(r,z^2)\frac{u_{m_\ell}^{2^*-1}\xi_{\ell}}{|y|}\nu_{k+j},
\end{equation}
and
\begin{equation}\label{po2}
\int_{\Omega_1}\frac{u_{m_\ell}^{2^*-1} \xi_{\ell} }{|y|}\langle\nabla K,x-(0,\zeta_1)\rangle=-\langle\nu,(0,\zeta_1)\rangle \Big(\int_{\partial\Omega_1} \frac{K(r,z^2)}{|y|}u_{m_\ell}^{2^*-1}\xi_{\ell}  -\int_{\partial\Omega_1} \nabla u_{m_\ell}\cdot\nabla\xi_{\ell}\Big).
\end{equation}
Combining \eqref{po1} and \eqref{po2}, we obtain
\begin{equation}\label{poid1}
\int_{\Omega_1}\frac{\partial K(r,z^2)}{\partial z_j}\frac{u_{m_\ell}^{2^*-1}\xi_{\ell}}{|y|}
=-\frac{ \nu_{k+j} }{\langle\nu,(0,\zeta_1)\rangle} \int_{\Omega_1}\frac{u_{m_\ell}^{2^*-1} \xi_{\ell} }{|y|}\langle\nabla K,x-(0,\zeta_1)\rangle.
\end{equation}

Next, we give the estimate to the terms of both side of \eqref{poid1}. By symmetry, we have
$$
\int_{\R^n}\frac{U^{2^*-1}}{|y|}\Phi_{i}=  \int_{\R^n} -\Delta U \Phi_{i}= \int_{\R^n} -\Delta\Phi_{i}U  =(2^*-1) \int_{\R^n}\frac{U^{2^*-1}}{|y|}\Phi_{i}=0,
$$
for $i=1,3,4,\cdots,n-k$. Then
\begin{align*}
    \int_{\Omega_1} \frac{u_{m_\ell}^{2^*-1}\xi_{\ell}}{|y|}=&\int_{\R^n} \frac{U^{2^*-1}}{|y|}\Big( b_{0,\ell}\Phi_{0} + b_{1,\ell}\Phi_{1}+\sum_{i=3}^{n-k}b_{i,\ell}\Phi_{i}+\mu_{m_\ell}^{-\frac{n-2}{2}} \xi_\ell^*(\mu_{m_\ell}^{-1}x+(0,\zeta_1))\Big)+O\Big(\frac{1}{\mu_{m_\ell}^2}\Big)\\
    =&b_{0,\ell}\int_{\R^n} \frac{U^{2^*-1}}{|y|}\Phi_{0}+O\Big(\frac{1}{\mu_{m_\ell}^{\frac{n-2}{2}-\tau}}\Big)+O\Big(\frac{1}{\mu_{m_\ell}^2}\Big)\\
    =&O\Big(\frac{1}{\mu_{m_\ell}^{\min\{\frac{n-2}{2}-\tau,2\}}}\Big),
\end{align*}
and
\begin{align*}
    \int_{\Omega_1} u_{m_\ell}\xi_{\ell}=&\frac{1}{\mu_{m_\ell}^2}\int_{\R^n}  U\Big( b_{0,\ell}\Phi_{0} + b_{1,\ell}\Phi_{1}+\sum_{i=3}^{n-k}b_{i,\ell}\Phi_{i}+\mu_{m_\ell}^{-\frac{n-2}{2}} \xi_\ell^*(\mu_{m_\ell}^{-1}x+(0,\zeta_1))\Big)+O\Big(\frac{1}{\mu_{m_\ell}^{2+2\tau}}\Big)\\
    =&\frac{b_{0,\ell}}{\mu_{m_\ell}^2}\int_{\R^n} U\Phi_{0}+O\Big(\frac{1}{\mu_{m_\ell}^{2+2\tau}}\Big).
\end{align*}
Since $\nabla K (0,\zeta_1)=O(|(\bar r_{m_\ell},\bar z^2_{m_\ell})-(r_0,z^2_0)|),$ then
\begin{equation}\label{po11}
\begin{aligned}
      & \int_{\Omega_1}\frac{\partial K(r,z^2)}{\partial z_j}\frac{u_{m_\ell}^{2^*-1}\xi_{\ell}}{|y|}\\
      =& \int_{\Omega_1}\Big(\frac{\partial K(r,z^2)}{\partial z_j}-\frac{\partial K(0,\zeta_1)}{\partial z_j}\Big)\frac{u_{m_\ell}^{2^*-1}\xi_{\ell}}{|y|}+O\Big(\frac{1}{\mu_{m_\ell}^{\min\{\frac{n}{2}-\bar\theta-\tau,3-\bar\theta\}}}\Big)\\
      =&\sum_{i=1,i\neq 2}^{n-k}\frac{b_{i,\ell}}{\mu_{m_\ell}}\frac{\partial^2K(0,\zeta_1)}{\partial z_j\partial z_i}\int_{\R^n}\frac{U^{2^*-1}\Phi_i}{|y|}z_i+\frac{b_{0,\ell}}{2(n-k)\mu_{m_\ell}^2}\frac{\partial \Delta K(0,\zeta_1)}{\partial z_j}\int_{\R^n}\frac{U^{2^*-1}\Phi_0}{|y|}|z|^2+O\Big(\frac{1}{\mu_{m_\ell}^{2+\sigma}}\Big).
\end{aligned}
\end{equation}
On the other hand, we have
\begin{equation}\label{po12}
\begin{aligned}
      &\int_{\Omega_1}\frac{u_{m_\ell}^{2^*-1} \xi_{\ell} }{|y|}\langle\nabla K,x-(0,\zeta_1)\rangle\\
      =&\int_{\Omega_1}\frac{u_{m_\ell}^{2^*-1} \xi_{\ell} }{|y|}\langle\nabla K(r,z^2)-\nabla K(0,\zeta_1),x-(0,\zeta_1)\rangle+O\Big(\frac{1}{\mu_{m_\ell}^{2+\sigma}}\Big)\\
    =&\frac{b_{0,\ell} \Delta K(0,\zeta_1)}{(n-k)\mu_{m_\ell}^2} \int_{\R^n}\frac{U^{2^*-1}\Phi_0}{|y|}|z|^2+O\Big(\frac{1}{\mu_{m_\ell}^{2+\sigma}}\Big).
\end{aligned}
\end{equation}
Therefore, from \eqref{poid1}--\eqref{po12}, we can obtain that for $j=1,3,4,\cdots,n-k$,
\begin{align}\label{bin}
    & \sum_{i=1,i\neq 2}^{n-k}b_{i,\ell}\frac{\partial^2K(0,\zeta_1)}{\partial z_j\partial z_i}\int_{\R^n}\frac{U^{2^*-1}\Phi_i}{|y|}z_i\nonumber\\
    =& -b_{0,\ell}\Big(\Big(\frac{ \nu_{k+j} }{\langle\nu,(0,\zeta_1)\rangle}\Delta K(0,\zeta_1)+\frac{1}{2}\frac{\partial \Delta K(0,\zeta_1)}{\partial z_j}\Big)\frac{ 1}{(n-k)\mu_{m_\ell}}\int_{\R^n}\frac{U^{2^*-1}\Phi_0}{|y|}|z|^2\Big)
    +O\Big(\frac{1}{\mu_{m_\ell}^{1+\sigma}}\Big).
\end{align}
The assumption ${\bf (K_3)}$ indicates that linear system \eqref{bin} is solvable, hence
by the boundness of $b_{0,\ell}$, we know that $b_{i,\ell}=O(\frac{1}{\mu_{m_\ell}})=o(1)$, $i=1,3,4,\cdots,n-k$.

{\bf Step 2.} We claim that $b_{0,\ell}\to 0$. In order to get an estimate to $b_{0,\ell}$, we apply the local Pohozaev identity \eqref{pohozeav2} in $B_{{\delta}/{{m_\ell}}}(0,\zeta_1)$, where $\delta>0$ is a fixed small constant, we have
\begin{equation}\label{poho1}
\begin{aligned}
     &\int_{B_{{\delta}/{{m_\ell}}}(0,\zeta_1)}\frac{u_{m_\ell}^{2^*-1} \xi_\ell }{|y|}\langle\nabla K,x-(0,\zeta_1)\rangle
     \\=&\int_{\partial B_{{\delta}/{{m_\ell}}}(0,\zeta_1)} \frac{K(r,z^2)}{|y|}u_{m_\ell}^{2^*-1}\xi_\ell \langle\nu,x-(0,\zeta_1)\rangle
\\&+\int_{\partial B_{{\delta}/{{m_\ell}}}(0,\zeta_1)} \Big(\frac{\partial u_{m_\ell}}{\partial \nu}\langle\nabla \xi_\ell,x-(0,\zeta_1)\rangle+\frac{\partial \xi_\ell}{\partial \nu}\langle\nabla u_{m_\ell},x-(0,\zeta_1)\rangle\Big)
\\&-\int_{\partial B_{{\delta}/{{m_\ell}}}(0,\zeta_1)} \nabla u_{m_\ell} \cdot\nabla\xi_\ell\langle\nu,x-(0,\zeta_1)\rangle+\frac{n-2}{2}\int_{\partial B_{{\delta}/{{m_\ell}}}(0,\zeta_1)}\Big(u_{m_\ell}\frac{\partial \xi_\ell}{\partial \nu}+\xi_\ell \frac{\partial u_{m_\ell}}{\partial \nu}\Big).
\end{aligned}
\end{equation}
A direct computation shows that
\begin{equation}\label{pohole}
\begin{aligned}
    &\int_{B_{{\delta}/{{m_\ell}}}(0,\zeta_1)}\frac{u_{m_\ell}^{2^*-1} \xi_\ell }{|y|}\langle\nabla K,x-(0,\zeta_1)\rangle\\
    =&\frac{b_{0,\ell} \Delta K(0,\zeta_1)}{(n-k)\mu_{m_\ell}^2} \int_{\R^n}\frac{U^{2^*-1}\Phi_0}{|y|}|z|^2+O\Big(\frac{1}{\mu_{m_\ell}^{2+\sigma}}\Big)\\
    =&-\frac{4b_{0,\ell}}{n-2}\frac{ \Delta K(0,\zeta_1)}{2^*(n-k)\mu_{m_\ell}^2} \int_{\R^n}\frac{U^{2^*}}{|y|}|z|^2+O\Big(\frac{1}{\mu_{m_\ell}^{2+\sigma}}\Big),
\end{aligned}
\end{equation}
and
$$
\begin{aligned}
    &\Big|\int_{\partial B_{{\delta}/{{m_\ell}}}(0,\zeta_1)} \frac{K(r,z^2)}{|y|}u^{2^*-1}\xi \langle\nu,x-(0,\zeta_1)\rangle\Big|\\
    \le &C\int_{s^2+t^2=(\frac{\mu_{m_\ell}\delta }{{m_\ell}})^2} \frac{s^{k-2}t^{n-k-1}}{(1+s+t)^{2n-3}}+\int_{s^2+t^2=(\frac{\mu_{m_\ell}\delta }{{m_\ell}})^2} \frac{1}{\mu_{m_\ell}^2}\frac{s^{k-2}t^{n-k-1}}{(1+s+t)^{n-1}}\\
    =&O\Big(\frac{m_\ell^n}{\mu_{m_\ell}^{n}}\Big)+O\Big(\frac{1}{\mu_{m_\ell}^{2+\sigma}}\Big)=O\Big(\frac{1}{\mu_{m_\ell}^{2+\sigma}}\Big).
\end{aligned}
$$
Define
\begin{equation}\label{J}
\begin{aligned}
    J(u,\xi,d)=&\int_{\partial B_{d}(0,\zeta_1)} \Big(\frac{\partial u}{\partial \nu}\langle\nabla \xi,x-(0,\zeta_1)\rangle+\frac{\partial \xi}{\partial \nu}\langle\nabla u,x-(0,\zeta_1)\rangle\Big)\\
&-\int_{\partial B_{d}(0,\zeta_1)} \nabla u \cdot\nabla\xi\langle\nu,x-(0,\zeta_1)\rangle+\frac{n-2}{2}\int_{\partial B_{d}(0,\zeta_1)}\Big(u\frac{\partial \xi}{\partial \nu}+\xi \frac{\partial u}{\partial \nu}\Big).
\end{aligned}
\end{equation}

Denote that $G(\tilde x, x)=((n-2)\omega_n)^{-1}|\tilde x-x|^{2-n}$ be the Green's function of the operator $-\Delta$  in $\mathbb R^n$, where $\omega_n$ is the  volume of unit ball in $\mathbb{R}^{n}$. Let
$$
\partial_j G(\tilde x,x)=\frac{\partial G(\tilde x,x)}{\partial \tilde z_j},\qquad\nabla_i G(\tilde x,x)=\frac{\partial G(\tilde x,x)}{\partial  z_i}.
$$

Then for any $0<\varepsilon<d<{\delta}/{m_\ell},$ we have
\begin{equation}\label{Jg1}
\begin{aligned}
    &J(G(\tilde x,(0,\zeta_1)),G(\tilde x,(0,\zeta_1)),d)\\
    =&2\int_{B_d(0,\zeta_1)}\Delta G(\tilde x,(0,\zeta_1))\langle\nabla G(\tilde x,(0,\zeta_1)), \tilde x-(0,\zeta_1)\rangle\\
    &+(n-2)\int_{B_d(0,\zeta_1)}\Delta G(\tilde x,(0,\zeta_1)) G(\tilde x,(0,\zeta_1))\\
    =&0.
\end{aligned}
\end{equation}
And,
\begin{equation}\label{Jgij}
\begin{aligned}
    &J(G(\tilde x,(0,\zeta_1)),G(\tilde x,(0,\zeta_j)),d)-J(G(\tilde x,(0,\zeta_1)),G(\tilde x,(0,\zeta_j)),\varepsilon)\\
    =&\int_{B_d(0,\zeta_1)\setminus B_\varepsilon(0,\zeta_1)}\Delta G(\tilde x,(0,\zeta_1))\langle\nabla G(\tilde x,(0,\zeta_j)), \tilde x-(0,\zeta_1)\rangle\\
    &+\int_{B_d(0,\zeta_1)\setminus B_\varepsilon(0,\zeta_1)}\Delta G(\tilde x,(0,\zeta_j))\langle\nabla G(\tilde x,(0,\zeta_1)), \tilde x-(0,\zeta_1)\rangle\\
    &+\frac{n-2}{2}\int_{B_d(0,\zeta_1)\setminus B_\varepsilon(0,\zeta_1)}\Delta G(\tilde x,(0,\zeta_1)) G(\tilde x,(0,\zeta_j)\\
    &+\frac{n-2}{2}\int_{B_d(0,\zeta_1)\setminus B_\varepsilon(0,\zeta_1)}\Delta G(\tilde x,(0,\zeta_j)) G(\tilde x,(0,\zeta_1))\\
    =&0.
\end{aligned}
\end{equation}
Thus, for $j=2,\cdots,{m_\ell},$
\begin{equation}\label{JGG}
\begin{aligned}
    J(G(\tilde x,(0,\zeta_1)),G(\tilde x,(0,\zeta_j)),d)
    =&\lim_{\varepsilon\to 0}J(G(\tilde x,(0,\zeta_1)),G(\tilde x,(0,\zeta_j)),\varepsilon)\\
    =&-\frac{n-2}{2}G((0,\zeta_1),(0,\zeta_j)).
\end{aligned}
\end{equation}
Similarly, we have for $i=1,\cdots,n-k,$
\begin{equation}\label{Jg11}
J(G(\tilde x,(0,\zeta_1)),\nabla_i G(\tilde x,(0,\zeta_1)),d)=0,
\end{equation}
and for $j=2,\cdots,{m_\ell},$
\begin{equation}\label{Jgij1}
\begin{aligned}
    J(G(\tilde x,(0,\zeta_1)),\nabla_i G(\tilde x,(0,\zeta_j)),d)
    =&\lim_{\varepsilon\to 0}J(G(\tilde x,(0,\zeta_1)),\nabla_i G(\tilde x,(0,\zeta_j)),\varepsilon)\nonumber\\
    =&-\frac{n-2}{2}\nabla_i G((0,\zeta_1),(0,\zeta_j)).
\end{aligned}
\end{equation}
It follows  from \eqref{J} that
\begin{equation}
    \frac{J(G(\tilde x,(0,\zeta_1)),\nabla_i G(\tilde x,(0,\zeta_j)),d)}{\mu_{m_\ell}}=o\Big(\frac{1}{\mu_{m_\ell}}\Big).
\end{equation}
Therefore, combining \eqref{J}--\eqref{Jgij1}, and using the result of Lemma 3.2 in \cite{LiuWangCylindrical2023}, we get
\begin{equation}\label{pohori}
\begin{aligned}
    &J\Big(u_{m_\ell},\xi_\ell,\frac{\delta}{{m_\ell}}\Big)\\
    =&2\frac{b_{0,\ell}A_1A_2}{\mu_{m_\ell}^{n-2}}J\Big( G(\tilde x,(0,\zeta_1)),\sum_{j=2}^{m_\ell} G(\tilde x,(0,\zeta_j)),\frac{\delta}{{m_\ell}} \Big)\\
    &+2\frac{b_{1,\ell}A_1A_3}{\mu_{m_\ell}^{n-1}}J\Big( G(\tilde x,(0,\zeta_1)),\sum_{j=2}^{m_\ell}(\cos \theta_j \nabla_1 G(\tilde x,(0,\zeta_j))+\sin \theta_j \nabla_2 G(\tilde x,(0,\zeta_j))),\frac{\delta}{{m_\ell}} \Big)\\
     &+2\sum_{i=3}^{n-k}\frac{b_{i,\ell}A_1A_3}{\mu_{m_\ell}^{n-1}}J\Big( G(\tilde x,(0,\zeta_1)),\sum_{j=2}^{m_\ell}\nabla_i G(\tilde x,(0,\zeta_j)),\frac{\delta}{{m_\ell}} \Big)\\
     &+2\sum_{i=1,i\neq 2}^{n-k}\frac{b_{i,\ell}A_1A_3}{\mu_{m_\ell}^{n-1}}J\Big( \sum_{j=2}^{m_\ell}G(\tilde x,(0,\zeta_j)),\nabla_i G(\tilde x,(0,\zeta_1)),\frac{\delta}{{m_\ell}} \Big)+o\Big(\frac{1}{\mu_{m_\ell}^2}\Big)\\
    =&-(n-2)\frac{b_{0,\ell}A_1A_2}{\mu_{m_\ell}^{n-2}}\sum_{j=2}^{m_\ell} G((0,\zeta_1),(0,\zeta_j))+o\Big(\frac{1}{\mu_{m_\ell}^2}\Big)\\
    =&-\frac{b_{0,\ell}A_1A_2}{\omega_n\mu_{m_\ell}^{n-2}}\sum_{j=2}^{m_\ell}\frac{1}{|\zeta_j-\zeta_1|^{n-2}}+o\Big(\frac{1}{\mu_{m_\ell}^{2}}\Big)\\
    =&-b_{0,\ell}\frac{A_1 A_4}{(2^*-1)\omega_n\mu_{m_\ell}^{2}}+o\Big(\frac{1}{\mu_{m_\ell}^{2}}\Big),
\end{aligned}
\end{equation}
where $A_1, A_2, A_3, A_4$ are the constants defined by
$$
\begin{aligned}
    A_1:=&\int_{\R^n}\frac{U^{2^*-1}}{|y|}>0, \\
    A_2:=&(2^*-1)\int_{\R^n}\frac{U^{2^*-2}\Phi_0}{|y|}<0,\\
    A_3:=&(2^*-1)\displaystyle\int_{\R^n}\frac{U^{2^*-1}\Phi_1 z_1}{|y|}>0, \\
    A_4:=&-\frac{\Delta K(r_0,z_0^2)}{2^*(n-k)}\int_{\R^n}\frac{|z|^2U^{2^*}}{|y|}>0,
\end{aligned}
$$
and the last equality comes from
$$
\int_{\R^n}(2^*-1)\frac{U_{\zeta_1,\mu_{m_\ell}}^{2^*-2}}{|y|}
\sum_{j=2}^{m_\ell} U_{\zeta_j,\mu_{m_\ell}}  \frac{\partial U_{\zeta_1,\mu_{m_\ell}}}{\partial \mu_{m_\ell}}=-\sum_{i=2}^{m_\ell}\frac{(2^*-1 )A_2}{\mu_{m_\ell}^{n-1}|\zeta_1-\zeta_j|^{n-2}}=-\frac{A_4}{\mu_{m_\ell}^3}.
$$

Combining \eqref{poho1}, \eqref{pohole} and \eqref{pohori} we have
$$
\Big(\frac{4}{n-2}+\frac{A_1}{(2^*-1)\omega_n}\Big)A_4\frac{b_{0,\ell}}{\mu_{m_\ell}^2}
=o\Big(\frac{1}{\mu_{m_\ell}^2}\Big).
$$
Thus we deduce that $b_{0,\ell}\to 0.$
\end{proof}

Finally, we will complete the proof of Proposition \ref{pro2.4} by giving the expansion of $u_{m_\ell}, \xi_\ell$ and their partial derivatives $\frac{\partial u_{m_\ell}}{\partial z_i}$, $\frac{\partial \xi_\ell}{\partial z_i}$ on $\partial B_{{\delta}/{{m_\ell}}}(0,\zeta_1)$. We have the following lemma.

\begin{lemma}\label{lemma2.5}
For a small constant $\delta>0$ fixed, we have for any $\tilde x \in \partial B_{{\delta}/{{m_\ell}}}(0,\zeta_1)$,
\begin{equation}\label{umn}
    u_{m_\ell}(\tilde x)=\frac{A_1}{\mu_{m_\ell}^{\frac{n-2}{2}}}\sum_{j=1}^{m_\ell} G(\tilde x,(0,\zeta_j))+O\Big(\frac{m_\ell^{n-2}}{\mu_{m_\ell}^{\frac{n-2}{2}+\frac{2}{n-2}}}\Big),
\end{equation}
\begin{equation}\label{pumn}
\frac{\partial u_{m_\ell}}{\partial \tilde z_l}(\tilde x) =\frac{A_1}{\mu_{m_\ell}^{\frac{n-2}{2}}}\sum_{j=1}^{m_\ell} \partial_l G(\tilde x,(0,\zeta_j))+O\Big(\frac{m_\ell^{n-3}}{\mu_{m_\ell}^{\frac{n-2}{2}+\frac{2}{n-2}}}\Big),
\end{equation}
and
\begin{equation}\label{xin}
\begin{aligned}
    &\xi_{\ell}(\tilde x)\\
    =&\frac{b_{0,\ell}A_2}{\mu_{m_\ell}^{\frac{n-2}{2}}}\sum_{j=1}^{m_\ell}G(\tilde x,(0,\zeta_j))
    +\frac{b_{1,\ell}A_3}{\mu_{m_\ell}^{\frac{n}{2}}}\sum_{j=1}^{m_\ell}( \cos \theta_j \nabla_1 G(\tilde x,(0,\zeta_j))+\sin \theta_j \nabla_2 G(\tilde x,(0,\zeta_j))\\
    &+\sum_{i=3}^{n-k}\frac{b_{i,\ell}A_3}{\mu_{m_\ell}^{\frac{n}{2}}}\sum_{j=1}^{m_\ell}  \nabla_i G(\tilde x,(0,\zeta_j))+O\Big(\frac{m_\ell^{n-2}}{\mu_{m_\ell}^{\frac{n-2}{2}+\frac{2}{n-2}}}\Big),
\end{aligned}
\end{equation}
\begin{equation}\label{pxin}
\begin{aligned}
   &\frac{\partial \xi_{\ell}}{\partial \tilde z_l}(\tilde x) \\
   =&\frac{b_{0,\ell}A_2}{\mu_{m_\ell}^{\frac{n-2}{2}}}\sum_{j=1}^{m_\ell}\partial_l G(\tilde x,(0,\zeta_j))+\frac{b_{1,\ell}A_3}{\mu_{m_\ell}^{\frac{n}{2}}}\sum_{j=1}^{m_\ell}\partial_l ( \cos \theta_j \nabla_1 G(\tilde x,(0,\zeta_j))+\sin \theta_j \nabla_2 G(\tilde x,(0,\zeta_j)))\\
    &+\sum_{i=3}^{n-k}\frac{b_{i,\ell}A_3}{\mu_{m_\ell}^{\frac{n}{2}}}\sum_{j=1}^{m_\ell}  \nabla_i \partial_l G(\tilde x,(0,\zeta_j))+O\Big(\frac{m_\ell^{n-3}}{\mu_{m_\ell}^{\frac{n-2}{2}+\frac{2}{n-2}}}\Big).
\end{aligned}
\end{equation}
\end{lemma}

\begin{proof}
Noting that
$$
u_{m_\ell}(\tilde x)=\int_{\R^n} G(\tilde x,x)K(r,z^2)\frac{u_{m_\ell}^{2^*-1}(x)}{|y|},
$$
and
\begin{equation}\label{pumnint}
\frac{\partial u_{m_\ell}}{\partial \tilde z_l}(\tilde x)=\int_{\R^n} \partial_l G(\tilde x,x)K(r,z^2)\frac{u_{m_\ell}^{2^*-1}(x)}{|y|}.
\end{equation}
With loss of generality, we assume that $y\in \Omega_1,$ and set $d={|\tilde x-(0,\zeta_1)|}/{2}$, we have
    \begin{align}
    u_{m_\ell}(\tilde x)=&\int_{B_{d}(\tilde x)\cup B_{d}(0,\zeta_1)}G(\tilde x,x)K(r,z^2)\frac{u_{m_\ell}^{2^*-1}}{|y|}\nonumber\\
    &+\int_{\Omega_1\setminus B_{d}(\tilde x)\setminus B_{d}(0,\zeta_1)}G(\tilde x,x)K(r,z^2)\frac{u_{m_\ell}^{2^*-1}}{|y|}\nonumber\\
    &+\sum_{j=2}^{m_\ell}\int_{ B_{{\delta}/{{m_\ell}}}(0,\zeta_j)}G(\tilde x,x)K(r,z^2)\frac{u_{m_\ell}^{2^*-1}}{|y|}\nonumber\\
    &+\sum_{j=2}^{m_\ell}\int_{\Omega_j\setminus B_{{\delta}/{{m_\ell}}}(0,\zeta_j)} G(\tilde x,x)K(r,z^2)\frac{u_{m_\ell}^{2^*-1}}{|y|}\nonumber\\
    :=&I_1+I_2+I_3+I_4.
\end{align}
We first compute $I_1$. Notice that
\begin{equation}\label{i11}
\begin{aligned}
       &\int_{B_{d}(\tilde x)}G(\tilde x,x)K(r,z^2)\frac{u_{m_\ell}^{2^*-1}}{|y|}\\
       \le&C\int_{B_{d}(\tilde x)}\frac{1}{|\tilde x-x|^{n-2}}\Big(\frac{\mu_{{m_\ell}}^{\frac{n+2}{2}}}{\mu_{{m_\ell}}|y|(1+\mu_{{m_\ell}}|y|+\mu_{{m_\ell}}|z-\zeta_1|)^{n}}\\
       &\quad\quad+\frac{\mu_{{m_\ell}}^{\frac{n+2}{2}}}{\mu_{{m_\ell}}|y|(1+\mu_{{m_\ell}}|y|+\mu_{{m_\ell}}|z-\zeta_1|)^{n-\sigma}}\frac{1}{\mu_{{m_\ell}}^\frac{2\sigma}{n-2}}\Big)\\
       \le&C\mu_{m_\ell}^{\frac{n-2}{2}}\frac{1}{(\mu_{m_\ell}d)^n}\int_{0}^{\mu_{m_\ell}d}\int_{0}^{2\pi}\frac{\cos^{k-2}\alpha\sin^{n-k-1}\alpha}{(\cos\alpha+\sin\alpha)^{n-2}}\, \text{d}\alpha\, \text{d}r\\
       \le& C\frac{1}{\mu_{m_\ell}^{2\frac{n-1}{n-2}-\frac{n-2}{2}}}.
\end{aligned}
\end{equation}
And  by Taylor expansion, for $x\in B_d(0,\zeta_j),$ we have
$$
G(\tilde x,x)=G(\tilde x,(0,\zeta_j))+\sum_{i=1}^{n-k}\nabla_i G(\tilde x,(0,\zeta_j))(x-(0,\zeta_j))_i+O\Big(\frac{|x-(0,\zeta_j)|^2}{|\tilde z-\zeta_j|^{n}}\Big).
$$
Then
\begin{align}\label{i12}
       &\int_{B_{d}(0,\zeta_1)}G(\tilde x,x)K(r,z^2)\frac{u_{m_\ell}^{2^*-1}}{|y|}\nonumber\\
       =&\int_{B_{d}(0,\zeta_1)}\Big(G(\tilde x,(0,\zeta_1))+\sum_{i=1}^{n-k}\nabla_i G(\tilde x,(0,\zeta_1))(x-(0,\zeta_1))_i \Big)K(r,z^2)\frac{u_{m_\ell}^{2^*-1}}{|y|}\nonumber\\
       &+O\Big(\int_{B_{d}(0,\zeta_1)}\frac{|x-(0,\zeta_1)|^2}{d^n}\frac{u_{m_\ell}^{2^*-1}}{|y|}\Big)\nonumber\\
       =&G(\tilde x,(0,\zeta_1))\frac{1}{\mu_{m_\ell}^{\frac{n-2}{2}}}\int_{\R^n}\frac{K(0,\zeta_1)U^{2^*-1}}{|y|}+O\Big(\frac{G(\tilde x,(0,\zeta_1))}{\mu_{m_\ell}^{\frac{n+2}{2}}}\Big)\nonumber\\
       &+O\Big(G(\tilde x,(0,\zeta_1))+ \sum_{i=1}^{n-k}\frac{|\nabla_i G(\tilde x,(0,\zeta_1))|}{m_\ell}\Big)\times \Big(\frac{1}{\mu_{m_\ell}^{\frac{n-2}{2}}}\int_{\mu_{m_\ell}d}^{+\infty}\int_{0}^{2\pi}\frac{r^{n-2}\cos\alpha^{k-2}\sin \alpha^{n-k-1}}{(1+r\cos\alpha+r\sin\alpha)^n}\,\text{d}\alpha\,\text{d}r\Big) \nonumber\\
       &+O\Big(\frac{\ln (d\mu_{m_\ell})}{d^n\mu_{m_\ell}^{\frac{n+2}{2}}}\Big)\nonumber\\
       =&G(\tilde x,(0,\zeta_1))\frac{1}{\mu_{m_\ell}^{\frac{n-2}{2}}}\int_{\R^n}\frac{K(0,\zeta_1)U^{2^*-1}}{|y|}+O\Big(\frac{G(\tilde x,(0,\zeta_1))+ \sum_{i=1}^{n-k}\frac{|\nabla_i G(\tilde x,(0,\zeta_1))|}{m_\ell}}{\mu_{m_\ell}^{\frac{n-2}{2}+\frac{2}{n-2}}}+\frac{1}{\mu_{m_\ell}^{\frac{2n}{n-2}-\frac{n-2}{2}}}\Big)\nonumber\\
       =&G(\tilde x,(0,\zeta_1))\frac{1}{\mu_{m_\ell}^{\frac{n-2}{2}}}\int_{\R^n}\frac{U^{2^*-1}}{|y|}+O\Big(\frac{m_\ell^{n-2}}{\mu_{m_\ell}^{\frac{n-2}{2}+\frac{2}{n-2}}}\Big).
\end{align}
   Therefore, combining \eqref{i11}--\eqref{i12}, we have
\begin{equation}\label{i1}
    I_1=\frac{G(\tilde x,(0,\zeta_1))}{\mu_{m_\ell}^{\frac{n-2}{2}}}\int_{\R^n}\frac{U^{2^*-1}}{|y|}+O\Big(\frac{m_\ell^{n-2}}{\mu_{m_\ell}^{\frac{n-2}{2}+\frac{2}{n-2}}}\Big).
\end{equation}
Similarly, we can also have
\begin{equation}\label{i2}
     I_2=\frac{\sum_{j=2}^{m_\ell} G(\tilde x,(0,\zeta_j))}{\mu_{m_\ell}^{\frac{n-2}{2}}}\int_{\R^n}\frac{U^{2^*-1}}{|y|}+O\Big(\frac{m_\ell^{n-2}}{\mu_{m_\ell}^{\frac{n-2}{2}+\frac{2}{n-2}}}\Big).
\end{equation}
For $I_3$ and $I_4$, we can calculate that
\begin{equation}\label{i3}
    |I_3|= O\Big(\frac{m_\ell^{n-2}}{\mu_{m_\ell}^{\frac{n-2}{2}+\frac{2}{n-2}}}\Big),
\end{equation}
and
\begin{equation}\label{i4}
    |I_4|= O\Big(\frac{m_\ell^{n-2}}{\mu_{m_\ell}^{\frac{n-2}{2}+\frac{2}{n-2}}}\Big).
\end{equation}
Therefore, combining \eqref{i1}--\eqref{i4}, we have \eqref{umn}. Similarly,  from \eqref{pumnint}, we can get \eqref{pumn}.

On the other hand, noting that
$$
\xi_{\ell}(\tilde x)=\int_{\R^n}(2^*-1) G(\tilde x,x)K(r,z^2)\frac{u_{m_\ell}^{2^*-2}(x)\xi_{\ell}( x)}{|y|},
$$
and
$$
\frac{\partial \xi_{\ell}}{\partial \tilde z_l}(\tilde x)=\int_{\R^n} (2^*-1)\partial_l G(\tilde x,x)K(r,z^2)\frac{u_{m_\ell}^{2^*-2}(x)\xi_{\ell}( x)}{|y|}.
$$
Then
\begin{equation}\label{j}
\begin{aligned}
    \xi_{\ell}(\tilde x)=&\int_{B_{d}(\tilde x)\cup B_{d}(0,\zeta_1)}(2^*-1)G(\tilde x,x)K(r,z^2)\frac{u_{m_\ell}^{2^*-2}\xi_{\ell}}{|y|}\\
    &+\int_{\Omega_1\setminus B_{d}(\tilde x)\setminus B_{d}(0,\zeta_1)}(2^*-1)G(\tilde x,x)K(r,z^2)\frac{u_{m_\ell}^{2^*-2}\xi_{\ell}}{|y|}\\
    &+\sum_{j=2}^{m_\ell}\int_{ B_{{\delta}/{{m_\ell}}}(0,\zeta_j)}(2^*-1)G(\tilde x,x)K(r,z^2)\frac{u_{m_\ell}^{2^*-2}\xi_{\ell}}{|y|}\\
    &+\sum_{j=2}^{m_\ell}\int_{\Omega_j\setminus B_{{\delta}/{{m_\ell}}}(0,\zeta_j)}(2^*-1)G(\tilde x,x)K(r,z^2)\frac{u_{m_\ell}^{2^*-2}\xi_{\ell}}{|y|}\\
    :=&J_1+J_2+J_3+J_4.
\end{aligned}
\end{equation}
Similar to the calculation of $u_{m_\ell},$ we can obtain
\begin{equation}\label{j1}
\begin{aligned}
     J_1=&\frac{b_{0,\ell}(2^*-1)}{\mu_{m_\ell}^{\frac{n-2}{2}}}G(\tilde x,(0,\zeta_1))\int_{\R^n}\frac{U^{2^*-2}\Phi_0}{|y|}\\
     &+\frac{b_{1,\ell}(2^*-1)}{\mu_{m_\ell}^{\frac{n}{2}}} \nabla_1 G(\tilde x,(0,\zeta_1))\int_{\R^n}\frac{U^{2^*-2}\Phi_1 z_1}{|y|}\\
    &+\sum_{i=3}^{n-k}\frac{b_{i,\ell}(2^*-1)}{\mu_{m_\ell}^{\frac{n}{2}}}\nabla_i G(\tilde x,(0,\zeta_1))\int_{\R^n}\frac{U^{2^*-2}\Phi_1 z_1}{|y|}+O\Big(\frac{m_\ell^{n-2}}{\mu_{m_\ell}^{\frac{n-2}{2}+\frac{2}{n-2}}}\Big),
\end{aligned}
\end{equation}
\begin{equation}\label{j2}
\begin{aligned}
      J_2=&\frac{b_{0,\ell}(2^*-1)}{\mu_{m_\ell}^{\frac{n-2}{2}}}\sum_{j=2}^{m_\ell}G(\tilde x,(0,\zeta_j))\int_{\R^n}\frac{U^{2^*-2}\Phi_0}{|y|}\\
      &+\frac{b_{1,\ell}(2^*-1)}{\mu_{m_\ell}^{\frac{n}{2}}}\sum_{j=2}^{m_\ell}( \cos \theta_j \nabla_1 G(\tilde x,(0,\zeta_j))+\sin \theta_j \nabla_2 G(\tilde x,(0,\zeta_j)))\int_{\R^n}\frac{U^{2^*-2}\Phi_1 z_1}{|y|}\\
    &+\sum_{i=3}^{n-k}\frac{b_{i,\ell}(2^*-1)}{\mu_{m_\ell}^{\frac{n}{2}}}\sum_{j=2}^{m_\ell}  \nabla_i G(\tilde x,(0,\zeta_j))\int_{\R^n}\frac{U^{2^*-2}\Phi_1 z_1}{|y|}+O\Big(\frac{m_\ell^{n-2}}{\mu_{m_\ell}^{\frac{n-2}{2}+\frac{2}{n-2}}}\Big),
\end{aligned}
\end{equation}
\begin{equation}\label{j34}
    |J_3|= O\Big(\frac{m_\ell^{n-2}}{\mu_{m_\ell}^{\frac{n-2}{2}+\frac{2}{n-2}}}\Big),\quad |J_4|= O\Big(\frac{m_\ell^{n-2}}{\mu_{m_\ell}^{\frac{n-2}{2}+\frac{2}{n-2}}}\Big).
\end{equation}
Combining \eqref{j}--\eqref{j34}, we have proved \eqref{xin}, and \eqref{pxin} can be proved similarly.
\end{proof}

Now, we give the proof of Theorem \ref{nondeg}.

\begin{proof}[Proof of Theorem \ref{nondeg}]
With the aid of the above lemmas and propositions, we are able to get a contradiction with $\|\xi_\ell\|_{*}=1.$ In fact, since
$$
|\xi_{\ell}(\tilde x)|
\le C\Big|\int_{\R^n} \frac{K(r,z^2)}{|\tilde x-x|^{n-2}}\frac{u_{m_\ell}^{2^*-2}(x)\xi_{\ell}( x)}{|y|}\Big|
\le C\|\xi_\ell\|_{*} \sum_{j=1}^{m_\ell}\frac{\mu_{m_\ell}^{\frac{n-2}{2}}}{(1 +\mu_{m_\ell}|\tilde y|+ \mu_{m_\ell}|\tilde z-{\zeta_j}|)^{\frac{n-2}{2}+\tau+\theta}},
$$
for some $\theta>0.$ Then we obtain
$$
{|\xi_{\ell}(\tilde x)|}\Big(\sum_{j=1}^{m_\ell}\frac{\mu_{m_\ell}^{\frac{n-2}{2}}}{(1 +\mu_{m_\ell}|\tilde y|+ \mu_{m_\ell}|\tilde z-{\zeta_j}| )^{\frac{n-2}{2}+\tau}}\Big)^{-1}
\le C\|\xi_\ell\|_{*} \frac{\sum_{j=1}^{m_\ell}\frac{\mu_{m_\ell}^{\frac{n-2}{2}}}{(1 +\mu_{m_\ell}|\tilde y|+ \mu_{m_\ell}|\tilde z-{\zeta_j}| )^{\frac{n-2}{2}+\tau+\theta}}}{\sum_{j=1}^{m_\ell}\frac{\mu_{m_\ell}^{\frac{n-2}{2}}}{(1 +\mu_{m_\ell}|\tilde y|+ \mu_{m_\ell}|\tilde z-{\zeta_j}| )^{\frac{n-2}{2}+\tau}}}.
$$
Noting that $\xi_\ell\to 0$ in $B_{{R}/{\mu_{m_\ell}}}(0,\zeta_1)$ and $\|\xi_\ell\|_{*}=1$, we know that
$$
{|\xi_{\ell}(\tilde x)|}\Big({\sum_{j=1}^{m_\ell}\frac{\mu_{m_\ell}^{\frac{n-2}{2}}}{(1 +\mu_{m_\ell}|\tilde y|+ \mu_{m_\ell}|\tilde z-{\zeta_j}|)^{\frac{n-2}{2}+\tau}}}\Big)^{-1}
$$
attains its maximum in $\R^n\setminus\bigcup_{j=1}^{m_\ell}B_{{R}/{\mu_{m_\ell}}}(0,\zeta_j)$. Thus,
$$
\|\xi_\ell\|_{*}\le o(1)\|\xi_\ell\|_{*}.
$$
So $\|\xi_\ell\|_{*}\to 0$ as $\ell\to +\infty$. This is a contradiction to $\|\xi_\ell\|_{*}=1$.
\end{proof}

\section{Construction of the new solutions}\label{constru}

In this section, we will construct a new kind of bubbling solutions. By Remark \ref{remark1}, we can obtain bubble solutions similar to \eqref{form}. In fact, let $m>0$ be a large even integer,
$$
	\zeta_{j}=\Big(\bar r\cos\frac{2(j-1)\pi}m, \bar r\sin\frac{2(j-1)\pi}m,0,0, \tilde z'\Big),\quad j=1,\cdots,m.
$$
Then under the condition $\bf{(K_1)}$-$\bf{(K_3)}$, we can prove that, in a similar way to the proof of Theorem A, there exist an integer $m_0>0$, such that for any even number $m>m_0$, problem \eqref{equation} has a solution $u_m$ of the form
$$
u_m=\overline{W}_{\bar r_m,  \tilde z'_m, \mu_m}+\phi_m,
$$
where $\overline{U}_{\zeta_j,\mu}=\eta U_{\zeta_j,\mu}, {W}_{\bar r, \tilde z',\mu}=\sum_{j=1}^{m}{U}_{\zeta_j,\mu}, \overline{W}_{\bar r, \tilde z',\mu}=\sum_{j=1}^{m}\overline{U}_{\zeta_j,\mu}$, $\phi_m \in {H_s}$, $(\bar r_m,  \tilde z'_m) \to (t_0, z'_0)$, $\mu_m \in [L_0 m^{\frac{n-2}{n-4} },L_1 m^{\frac{n-2}{n-4}}]$, and $\|\phi_m\|_{L^{\infty}(\R^n)}=o(\mu_m^{\frac{n-2}{2}})$.

Next, we will construct new cylindrical solution, by gluing bubbles at $(z_3,z_4)$-plane.
Let $q\ge m$ be a large integer.
Recall that
$$
p_{j}=\Bigl(0, 0,\bar t\cos\frac{2(j-1)\pi}q, \bar t\sin\frac{2(j-1)\pi}q,\bar z'\Bigr),\quad j=1,\cdots,q,
$$
where $ \bar z'\in \R^{n-k-4}$, $(\bar t, \bar z')\to (t_0, z_0')$. We introduce the weighted norms:
$$
\|u\|_{\tilde *}:= \sup_{x\in \R^n}\Big(\sum_{j=1}^{q}\frac{\lambda^{\frac{n-2}{2}}}{(1+\lambda |y|+\lambda |z-p_j|)^{\frac{n-2}{2}+\tau}}\Big)^{-1}|u(x)|,
$$
$$
\|f\|_{\widetilde {**}}:=\sup_{x\in \R^n}\Big(\sum_{j=1}^{q}\frac{\lambda^{\frac{n+2}{2}}}{\lambda |y|(1+\lambda |y|+\lambda |z-p_j|)^{\frac{n}{2}+\tau}}\Big)^{-1}|f(x)|,
$$
where $\tau=\frac{n-4}{n-2}.$
We aim to construct a solution of \eqref{equation} with the form
\begin{equation}
    v_q=u_m+\sum_{j=1}^{q} \eta U_{p_j,\lambda_q}+\psi_q,
\end{equation}
where $\psi_q\in X_s\cap D^{1,2}(\R^n)$ is a correction term, $X_s$ is defined as in \eqref{eq:Xs}. Throughout this section, we assume
\begin{equation}\label{para}
(\bar t, \bar z', \lambda)\in {\mathscr S_q}:= \Big\{(\bar t, \bar z', \lambda):|(\bar t, \bar z')-(t_0,z'_0)|\le\frac{1}{\lambda^{1-\bar \theta}},\, \lambda\in[L_0 q^{\frac{n-2}{n-4} },L_1 q^{\frac{n-2}{n-4}}]\Big\},
\end{equation}
with $\bar \theta \in (0,\frac{1-\epsilon_0}{2})$ and $\frac{n}{2}-\bar\theta-\tau>2$.

Consider the following linearized problem around $u_m+\sum_{j=1}^{q} \eta U_{p_j,\lambda_q}$:
\begin{equation}\label{lin}
\begin{cases}
-\Delta {\psi}-(2^*-1)K(t,z')\displaystyle\frac{(u_m+\overline{Y}_{\bar t,\bar z',\lambda})^{2^*-2}}{|y|}\psi=f+\displaystyle\sum\limits_{l=3}^{n-k}\displaystyle{c}_l\sum\limits_{j=1}^{k} \frac{\overline{Z}_{p_j,\lambda}^{2^*-2}  }{|y|}\overline{\mathbb{Z}}_{l j}\quad\text{ in }\,\R^n,
\\
\psi\in X_s, \quad \displaystyle\int_{\R^n} \frac{\overline{Z}_{p_j,\lambda}^{2^*-2} }{|y|} \overline{\mathbb{Z}}_{l j}\psi= 0,  \quad j= 1, \cdots, k, \quad l= 3,\cdots,n-k,
\end{cases}
\end{equation}
where $\overline{Z}_{p_j,\lambda}=\eta U_{p_j,\lambda}$, ${Y}_{\bar t,\bar z',\lambda}=\sum_{j=1}^{q}{Z}_{p_j,\lambda}$, $\overline{Y}_{\bar t,\bar z',\lambda}=\sum_{j=1}^{q}\overline{Z}_{p_j,\lambda}$, and the functions $\overline{\mathbb{Z}}_{l j}$ are defined as
$$
\overline{\mathbb{Z}}_{3j}=\frac{\partial \overline{Z}_{p_j,\lambda}}{\partial \lambda_q},
\quad
\overline{\mathbb{Z}}_{4j}=\frac{\partial \overline{Z}_{p_j,\lambda}}{\partial t_q},
\quad
\overline{\mathbb{Z}}_{l j}=\frac{\partial \overline{Z}_{p_j,\lambda}}{\partial z_j}, \, l=5,\cdots,n-k.
$$

\begin{lemma}\label{lem4.1}
Suppose that $n\ge 8$, $\frac{n+1}{2}\le k<n-3$, $K(x)$ satisfies ${\bf (K_1)}$-${\bf (K_3)}$,  and $(\bar t_q,\bar z_q', \lambda_q)\in {\mathscr S_q}$, $\psi_{q}$ solves \eqref{lin} for $f=f_{q}$. If  $\|f_{q}\|_{{\widetilde{**}}} \to 0$ as  $q \to +\infty$, then  $\|\psi_{q}\|_{{\tilde *}} \to 0$ as  $q \to +\infty$.
\end{lemma}

\begin{proof}
We argue  by contradiction. Suppose that there exists a sequence of ${\bar t_q}\to {t_0}$, ${\bar z'_q}\to{z'_0}$, ${\lambda_q}\in[\Lambda_2 q^\frac{n-2}{n-4},\Lambda_3 q^\frac{n-2}{n-4}]$, so that  $\psi_q$ solves \eqref{lin} with  $f=f_q$, $t= {\bar t_q}$, $z'={\bar z'_q}$, $\lambda= {\lambda_q}$, $\Vert f _{q}\Vert_{\widetilde {**}}\to 0 $, and  $\|\psi_q\|_{\tilde*}\ge C'>0$. Without loss of generality, we may assume that  $\|\psi_q\|_{\tilde*}=1 $.  We drop the subscript $q$ for simplicity.
Noting that
$$
L_m \psi=(2^*-1)\Big(\frac{(u_m+\overline{Y}_{\bar t,\bar z',\lambda})^{2^*-2}}{|y|}\psi-\frac{u_m^{2^*-2}}{|y|}\psi\Big)+f+\sum\limits_{l=3}^{n-k}{c}_l\sum\limits_{j=1}^{k} \frac{\overline{Z}_{p_j,\lambda}^{2^*-2}  }{|y|}\overline{\mathbb{Z}}_{l j}.
$$
Applying the Green's representation to $\psi$, we have
$$
\begin{aligned}
      |\psi(\tilde x)|
      \le& C\int_{\R^n}\frac{K(t,z')}{|\tilde x-x|^{n-2}} \frac{\overline{Y}_{\bar t,\bar z',\lambda}^{2^*-2}}{|y|}|\psi|+C\int_{\R^n}\frac1{|\tilde x-x|^{n-2}} |f|\\
      &+C\int_{\R^n}\frac1{|\tilde x-x|^{n-2}}\Big|\displaystyle\sum\limits_{l=3}^{n-k}\displaystyle{c}_l\sum\limits_{j=1}^{k} \frac{\overline{Z}_{p_j,\lambda}^{2^*-2}  }{|y|}\overline{\mathbb{Z}}_{l j}\Big|.
\end{aligned}
$$
Similar to the calculation in \cite{LiuWangCylindrical2023}, we can obtain
\begin{equation}\label{psi1}
\|{\psi}\|_{\tilde *} \le \Big(\|f\|_{\widetilde{**}} + \frac{ \sum_{j=1}^{q}(\frac{\lambda^{\frac{n-2}{2}}}{(1 + \lambda|\tilde y|+\lambda|\tilde z-{p_j}| )^{\frac{n-2}{2}+\tau+\theta}} )}{ \sum_{j=1}^{q}(\frac{\lambda^{\frac{n-2}{2}}}{(1 + \lambda|\tilde y|+\lambda|\tilde z-{p_j}| )^{\frac{n-2}{2}+\tau}} )}\|\psi\|_{\tilde *}+o_q(1)\Big),
\end{equation}
for some $\theta>0$ small enough.
Since  $\|{\psi}\|_{\tilde *}= 1 $, we obtain from \eqref{psi1} that there exist some positive constants  $R, \delta_1$ such that
\begin{equation}\label{bound1}
\|\lambda^{-\frac{n-2}{2}}\psi\|_{L^\infty(B_{{R}/{\lambda}}(0,p_{j}))}  \ge \delta_1> 0,
\end{equation}
for some  $j\in \{1,2, \cdots, q\}$. But $\tilde \psi(y):= \lambda^{-\frac{n-2}{2}}\psi(\lambda^{-1}y+(0,p_{j}))$  converges uniformly in any compact set to a solution $u$ of
\begin{equation}\label{ree1}
-\Delta v(x)- (2^*-1)\frac{U_{0,1}^{2^*-2}}{|y|} v(x)=0,\quad x=(y,z) \in \R^n,
\end{equation}
and $v$ is perpendicular to the kernel of \eqref{ree1}. As a result, $v=0$.  Together, with the non-degeneracy result, we deduce a contradiction to $\|\psi\|_{\tilde *}=1$.
\end{proof}

Now  we rewrite problem \eqref{lin} as the following perturbation problem:
\begin{equation}\label{linn}
\begin{cases}
\displaystyle {\bf L}_q \psi_q={\bf l}_q+{{\bf R}}(\psi_q)+\sum\limits_{l=3}^{n-k}{c}_l\sum\limits_{j=1}^{k} \frac{\overline{Z}_{p_j,\lambda_q}^{2^*-2}  }{|y|}\overline{\mathbb{Z}}_{l j}\quad \text{ in }\,\R^n,\\
\displaystyle\psi_q\in X_s, \quad \int_{\R^n} \frac{\overline{Z}_{p_j,\lambda_q}^{2^*-2} }{|y|} \overline{\mathbb{Z}}_{l j}\psi_q= 0,  \quad j= 1, \cdots, k, \, l= 3,\cdots,n-k,
\end{cases}
\end{equation}
where
$$
{\bf L}_q \psi_q:=-\Delta {\psi_q}-(2^*-1)K(t,z')\displaystyle\frac{(u_m+\overline{Y}_{\bar t_q,\bar z_q',\lambda_q})^{2^*-2}}{|y|}\psi_q,
$$
\begin{equation}\label{Lq}
\begin{aligned}
    {\bf l}_q:=&\frac{K(t,z')}{|y|}((u_m+\overline{Y}_{\bar t_q,\bar z_q',\lambda_q})^{2^*-1}-u_m^{2^*-1}-\overline{Y}_{\bar t_q,\bar z_q',\lambda_q}^{2^*-1})\\
    &+  \frac{K(t,z')\overline{Y}_{\bar t_q,\bar z_q',\lambda_q}^{2^*-1}-\sum\limits_{j=1}^{q}\eta U_{p_{j},\lambda_q}^{2^*-1}}{|y|}+\Delta\eta Y_{\bar t_q,\bar z_q',\lambda_q}+2\nabla\eta\nabla Y_{\bar t_q,\bar z_q',\lambda_q}\\
  :=&I_1+I_2+I_3+I_4,
\end{aligned}
\end{equation}
and
$$
{{\bf R}}(\psi_q):=\frac{K(t,z')}{|y|}( (u_m+\overline{Y}_{\bar t_q,\bar z_q',\lambda_q}+\psi_q)^{2^*-1}_{+}-(u_m+\overline{Y}_{\bar t_q,\bar z_q',\lambda_q})^{2^*-1}-(2^*-1)(u_m+\overline{Y}_{\bar t_q,\bar z_q',\lambda_q})^{2^*-2} \psi_q) .
$$

A standard argument leads to
\begin{lemma} \label{lem4.3}
Suppose that $n\ge 8$, $\frac{n+1}{2}\le k<n-3$, $K(x)$ satisfies ${\bf (K_1)}$-${\bf (K_3)},$ and $(\bar t_q,\bar z_q',\lambda_q)\in {\mathscr S_q}$, there exists $C>0$ such that
$$
\|{\bf R}(\psi_q)\|_{\widetilde{**}}\le C\|\psi_q\|_{\tilde *}^{2^*-1}.
$$
\end{lemma}

Next, we estimate $||{\bf l}_q||_{\widetilde{**}}.$

\begin{lemma} \label{lem2.4}
Suppose that $n\ge 8$, $\frac{n+1}{2}\le k<n-3$, $K(x)$ satisfies ${\bf (K_1)}$-${\bf (K_3)}$, $(\bar t_q,\bar z_q',\lambda_q)\in {\mathscr S_q}$, there exists $k_0>0 $ and $C>0$ such that for all $k \geq k_0$,
\begin{equation}\label{llk}
\|{\bf l_q}\|_{\widetilde{**}}\le \frac{C}{\lambda_q^{\frac{n-2\tau}{n-2}} }.
\end{equation}
\end{lemma}

\begin{proof}
Define
$$
\begin{aligned}
\widetilde\Omega_j:=\Big\{&x:x=(y,z^1,z_3,z_4,z')\in\R^k\times\R^2\times\R\times\R\times\R^{n-k-4},\\
&\Big\langle\frac{(z_3,z_4)}{|(z_3,z_4)|},\Big(\cos\frac{2(j-1)\pi}{q},\sin\frac{2(j-1)\pi}{q}\Big)\Big\rangle\ge\cos\frac{\pi}{q}\Big\}.
\end{aligned}
$$
We may assume $x\in \widetilde\Omega_{1}$ without loss of generality. Noting that for $x\in\widetilde\Omega_1\cap B_{\lambda_q^{-{1}/{2}}}(0,p_j),$ from Lemma \ref{lemma2.1}, we have
$$
|u_m|\le C \frac{mr_0^2}{\mu^{\frac{n-2}{2}}}\le C.
$$
Then
\begin{equation}\label{ii11}
\begin{aligned}
    |I_1|\le& C \frac{K(t,z')}{|y|}\Big(\Big(u_m+\sum_{j=2}^q \overline{U}_{p_j,\lambda_q}\Big)^{2^*-1}+\overline{U}_{p_1,\lambda_q}^{2^*-2} \Big(u_m+\sum_{j=2}^q \overline{U}_{p_j,\lambda_q}\Big)\Big) + \frac{C}{|y|}\\
    \le& C\frac{1}{|y|}\Big(\sum_{j=2}^q \overline{U}_{p_j,\lambda_q}\Big)^{2^*-1}+U_{p_1,\lambda_q}^{2^*-2} \Big(C+\sum_{j=2}^q \overline{U}_{p_j,\lambda_q}\Big)+\frac{C}{|y|}\\
    \le&C\Big(\Big(\frac{q}{\lambda_q}\Big)^{\frac{n}{2}-\tau} +\frac{1}{\lambda_q^{\frac{n}{4}-\frac{\tau}{2}}}\Big)\frac{\lambda_q^{\frac{n+2}{2}}}{\lambda_q|y|(1 +\lambda_q|y|+ \lambda_q|z-{p_1}|)^{\frac{n}{2}+\tau}}\\
           &+C\Big(\frac{q}{\lambda_q}\Big)^{(\frac{n-2}{2}-\frac{n-2}{n}\tau)\frac{n}{n-2}}    \Big(\sum_{j=2}^{q}\frac{\lambda_q^{\frac{n+2}{2}}}{\lambda_q|y|(1 +\lambda_q|y|+ \lambda_q|z-{p_j}| )^{\frac{n}{2}+\tau}}\Big).\\
\end{aligned}
\end{equation}
On the other hand, we consider the case that $x\notin \bigcup_{j=1}^q( \widetilde\Omega_1\cap B_{\lambda_q^{-{1}/{2}}}(0,p_j))$. Without loss of generality, we may assume $x\in\widetilde\Omega_1 \setminus B_{\lambda_q^{-{1}/{2}}}(0,p_j),$ then $U_{p_1,\lambda_q}\le C$.
Thus,
\begin{align}\label{ii12}
    |I_1|\le& C \frac{K(t,z')}{|y|}\Big(\Big(\sum_{j=1}^q \overline{U}_{p_j,\lambda_q}\Big)^{2^*-1}+ \sum_{j=1}^q \overline{U}_{p_j,\lambda_q}\Big)\nonumber\\
    \le&\frac{C}{|y|}\Big(\overline{U}_{p_1,\lambda_q}^{2^*-1}+ U_{p_1,\lambda_q}^{2^*-2}\Big(\sum_{j=2}^q \overline{U}_{p_j,\lambda_q}\Big) +  \sum_{j=1}^q \overline{U}_{p_j,\lambda_q}\Big)\nonumber\\
    \le& C\Big(\Big(\frac{q}{\lambda_q}\Big)^{(\frac{n-2}{2}-\frac{n-2}{n}\tau)\frac{n}{n-2}} +\frac{1}{\lambda_q^{\frac{n}{2}-2-\tau}}\Big) \Big(\sum_{j=2}^{q}\frac{\lambda_q^{\frac{n+2}{2}}}{\lambda_q|y|(1 +\lambda_q|y|+ \lambda_q|z-{p_j}| )^{\frac{n}{2}+\tau}}\Big)\nonumber\\
    &+C\frac{1}{\lambda_q^{\frac{n}{4}-\frac{\tau}{2}}}\frac{\lambda_q^{\frac{n+2}{2}}}{\lambda_q|y|(1 +\lambda_q|y|+ \lambda_q|z-{p_1}|)^{\frac{n}{2}+\tau}}.
\end{align}
Combining \eqref{ii11}--\eqref{ii12}, we have
\begin{equation}\label{3211}
\|I_1\|_{\widetilde{**}}\le C\max\Big\{ \Big(\frac{q}{\lambda_q}\Big)^{\frac{n}{2}-\tau},  \frac{1}{\lambda_q^{\frac{n-4}{2}-\tau}}, \frac{1}{\lambda_q^{\frac{n}{4}-\frac{\tau}{2}}}\Big\}\le \frac{C}{\lambda_q^{\frac{n-2\tau}{n-2}}}.
\end{equation}
Similar to the calculation in Lemma \ref{lemma2.3}, we have
\begin{equation}\label{3222}
\|I_2\|_{\widetilde{**}} \le  \frac{C}{\lambda_q^{\frac{n-2\tau}{n-2}}},
\end{equation}
and
\begin{equation}\label{3223}
\|I_3\|_{\widetilde{**}}+\|I_4\|_{\widetilde{**}} \le  \frac{C}{\lambda_q^{\frac{n-2}{2}-\tau}}.
\end{equation}
Hence, from \eqref{3211}--\eqref{3223}, we deduce that \eqref{llk} holds.
\end{proof}

By Fredholm alternative, and contraction mapping principle, we have the solvability theory for the  linearized  problem \eqref{linn} by a standard argument:

\begin{proposition}\label{pro4.5}
Suppose that $n\ge 8$, $\frac{n+1}{2}\le k<n-3$, $K(x)$ satisfies ${\bf (K_1)}$-${\bf (K_3)},$ and $(\bar t_q,\bar z_q',\lambda_q)\in {\mathscr S_q}$. There exists an integer  $q_0>0$ large enough, such that for each $q \ge q_0$, problem \eqref{linn} has a unique solution $\psi_q$ satisfying
\begin{equation} \label{phik}
\|\psi_q\|_{\tilde *} \le   \frac{C}{\lambda_q^{\frac{n-2}{2}-\tau}},\quad |c_{l}|\le  \frac{C}{\lambda_q^{\frac{n-2}{2}-\tau}}.
\end{equation}
\end{proposition}
	
Next, we have the following proposition which is necessary to choose proper $(\bar t, \bar z', \lambda)$ such that $u_m+\sum_{j=1}^{q} \overline{U}_{p_j,\lambda_q}+\psi_q$ be the solution of \eqref{equation}.

\begin{proposition}\label{pro3.5}
Suppose that  $n\ge 8$, $\frac{n+1}{2}\le k<n-3$, $K(x)$ satisfies ${\bf (K_1)}$-${\bf (K_3)}, $ and $(\bar t, \bar z', \lambda)$ satisfies
\begin{equation}\label{re1}
\int_{B_\rho}\Big(-\Delta v_q-K(t,z')\frac{(v_q)^{2^*-1}_+}{|y|}\Big)\langle x,\nabla v_q\rangle=0,
\end{equation}
\begin{equation}\label{re2}
\int_{B_\rho}\Big(-\Delta v_q-K(t,z')\frac{(v_q)^{2^*-1}_+}{|y|}\Big)\frac{\partial v_q}{\partial z_j}=0,\quad j=5,\cdots,n-k,
\end{equation}
and
\begin{equation}\label{re3}
\int_{\R^n}\Big(-\Delta v_q-K(t,z')\frac{(v_q)^{2^*-1}_+}{|y|}\Big)\frac{\partial \overline{Y}_{\bar t,\bar z',\lambda }}{\partial \lambda}=0,
\end{equation}
where $B_\rho:=\{(y,z^*,z')\in \R^{k}\times\R^4\times\R^{n-k-4}:|(|y|,|z^*|,z')-(0,t_0,z'_0)|\le\rho\}$ with $\rho \in (2\bar\delta, 5\bar\delta)$, $v_q=u_m+\sum_{j=1}^{q} U_{p_j,\lambda_q}+\psi_q$ is gotten from Proposition \ref{pro4.5}. Then
$$
c_l=0, \quad l=3,\cdots,n-k.
$$
\end{proposition}

\begin{proof}
Notice that
\begin{equation}\label{re4}
\int_{B_\rho}\Big(-\Delta u_m-K(t,z')\frac{(u_m)^{2^*-1}_+}{|y|}\Big)\langle x,\nabla u_m\rangle=0,
\end{equation}
and
$$
\psi_q\in X_s, \quad \int_{\R^n} \frac{\overline{Z}_{p_j,\lambda_q}^{2^*-2} }{|y|} \overline{\mathbb{Z}}_{l j}\psi_q= 0,  \quad j= 1, \cdots, k, \, l= 3,\cdots,n-k.
$$

Then \eqref{re1} is equivalent to
\begin{align}\label{equiv1}
&\int_{B_\rho}\Big(-\Delta u_q-K(t,z')\frac{(u_q)^{2^*-1}_+}{|y|}\Big)\langle x,\nabla u_q\rangle\nonumber\\
=&\int_{B_\rho}K(t,z')\frac{(v_q)^{2^*-1}_+-(u_m)^{2^*-1}_+ }{|y|}\langle x,\nabla u_q\rangle\nonumber\\
=&O\Big(\int_{B_\rho}K(t,z')\frac{u_m^{2^*-2}u_q+u_q^{2^*-1}}{|y|}\langle x,\nabla u_q\rangle\Big)=O(q)=o(q\lambda_q^2).
\end{align}
Similarly, \eqref{re2} is equivalent to
\begin{align}\label{equiv2}
&\int_{B_\rho}\Big(-\Delta u_q-K(t,z')\frac{(u_q)^{2^*-1}_+}{|y|}\Big)\frac{\partial v_q}{\partial z_j}\nonumber\\
=&\int_{B_\rho}K(t,z')\frac{(v_q)^{2^*-1}_+-(u_m)^{2^*-1}_+}{|y|}\frac{\partial v_q}{\partial z_j}=O(q)=o(q\lambda_q^2),
\end{align}
and \eqref{re3} is equivalent to
\begin{align}\label{equiv3}
&\int_{\R^n}\Big(-\Delta u_q-K(t,z')\frac{(u_q)^{2^*-1}_+}{|y|}\Big)\frac{\partial \overline{Y}_{\bar t,\bar z',\lambda}}{\partial \lambda}
=o\Big(\frac{q}{\lambda_q^{2}}\Big).
\end{align}

By similar argument of Proposition 3.1 in \cite{LiuWangCylindrical2023}, we can calculate from \eqref{equiv1}--\eqref{equiv3} that
\begin{equation}\label{cc1}
    c_4(a_3+o(1))=o\Big(\frac{1}{\lambda_q^2}\Big)c_3+\sum_{l=5}^{n-k}c_l(b_l+o(1)),
\end{equation}
and
\begin{equation}\label{cc2}
    c_j(a_4+o(1))=o\Big(\frac{1}{\lambda_q^2}\Big)c_3+o(1)\sum_{l=4,l\neq j}^{n-k} c_l,\quad j=5,\cdots,n-k,
\end{equation}
for some constants $a_3>0, a_4<0$, and $b_l\neq 0, l=5,\cdots,n-k$.
Then we deduce from \eqref{cc1}--\eqref{cc2} that
\begin{equation}\label{cc3}
    c_j=o\Big(\frac{1}{\lambda_q^2}\Big)c_3, \quad j=4,\cdots,n-k.
\end{equation}
On the other hand, we have from \eqref{equiv3} that $c_3$ satisfies that
\begin{equation}\label{cc4}
    \Big(a_5\frac{q}{\lambda_q^{2}}+ o\Big(\frac{q}{\lambda_q^{2}}\Big)\Big)c_3=o\Big(\frac{q}{\lambda_q^{2}}\Big),
\end{equation}
for $a_5>0.$ Thus, from \eqref{cc3}--\eqref{cc4}, we have $c_l=0, l=3,\cdots,n-k.$

\end{proof}
  For the construction of  new solutions, we can proceed exactly as in \cite{LiuWangCylindrical2023}. For readers convenience,  we give the sketch of the proof through the following lemmas and omit the detailed process.
 \begin{lemma}\label{lem4.6}
Suppose that  $n\ge 8$, $\frac{n+1}{2}\le k<n-3$, $K(x)$ satisfies ${\bf (K_1)}$-${\bf (K_3)}, $ then for $j=5,\cdots,n-k$, we have
\begin{equation}
    \int_{B_\rho}\frac{\partial K(t,z')}{\partial z_j}\frac{(v_q)_+^{2^*-1}}{|y|}=q\Big(\frac{\partial K(\bar t,\bar z')}{\partial \bar z_j}\int_{\R^N}\frac{U_{0,1}^{2^*}}{|y|}+o\Big(\frac{1}{q^{1/2}}\Big)\Big),
\end{equation}
and
\begin{equation}
    \int_{B_\rho}t\frac{\partial K(t,z')}{\partial t}\frac{(v_q)_+^{2^*-1}}{|y|}=q\Big(\bar t\frac{\partial K(\bar t,\bar z')}{\partial \bar t}\int_{\R^N}\frac{U_{0,1}^{2^*}}{|y|}+o\Big(\frac{1}{q^{1/2}}\Big)\Big).
\end{equation}
 \end{lemma}
\begin{lemma}\label{lem4.7}
    Suppose that  $n\ge 8$, $\frac{n+1}{2}\le k<n-3$, $K(x)$ satisfies ${\bf (K_1)}$-${\bf (K_3)}, $ then \eqref{re1} is equivalent to
    \begin{equation}\label{rre1}
        \int_{B_\rho}t\frac{\partial K(t,z')}{\partial t}\frac{(v_q)_+^{2^*-1}}{|y|}=o\Big(\frac{q}{\lambda_q^2}\Big).
    \end{equation}
    \eqref{re2} is equivalent to
    \begin{equation}\label{rre2}
        \int_{B_\rho}\frac{\partial K(t,z')}{\partial z_j}\frac{(v_q)_+^{2^*-1}}{|y|}=o\Big(\frac{q}{\lambda_q^2}\Big),\quad j=5,\cdots,n-k.
    \end{equation}
    And \eqref{re3} is equivalent to
    \begin{equation}\label{rre3}
        q\Big(\frac{C_1}{\lambda_q^3}-\frac{C_2 q^{n-2}}{\lambda_q^{n-1}}+o\Big(\frac{1}{\lambda_q^{3}}\Big)\Big)=0,
    \end{equation}
    where $C_1$ and $C_2$ are some positive constants.
\end{lemma}

 Define the energy functional:
$$
\mathcal{E}(v):=\frac{1}{2}\int_{\R^n}|\nabla v|^2-\frac{1}{2^*}\int_{\R^n}K(t,z')\frac{(v_+)^{2^*}}{|y|}.
$$
Now, we can give the proof of Theorem \ref{newsolution}.
\begin{proof}[Proof of Theorem \ref{newsolution}] We denote that
$$
F(\bar t_q,\bar z'_q,\lambda_q):=\mathcal{E}(u_m+\sum_{j=1}^{q} \overline{U}_{p_j,\lambda_q}+\psi_q),
$$
then by basic calculation, we have
\begin{align}
    &F(\bar t_q,\bar z'_q,\lambda_q)\nonumber\\
    =&\mathcal{E}\Big(u_m+\sum_{j=1}^{q} \overline{U}_{p_j,\lambda_q}\Big)+o\Big(\frac{q}{\lambda_q^2}\Big)\nonumber\\
    =&\mathcal{E}(u_m)+\mathcal{E}\Big(\sum_{j=1}^{q} \overline{U}_{p_j,\lambda_q}\Big)+o\Big(\frac{q}{\lambda_q^2}\Big)\nonumber\\
    &-\int_{\R^n}
    \frac{K(t,z')}{|y|}\Big(\Big(u_m+\sum_{j=1}^{q} \overline{U}_{p_j,\lambda_q}\Big)^{2^*}-u_m^{2^*}-\Big(\sum_{j=1}^{q} \overline{U}_{p_j,\lambda_q}\Big)^{2^*}-2^*\Big(\sum_{j=1}^{q}  u_m^{2^*-1}\overline{U}_{p_j,\lambda_q}\Big)\Big)\nonumber\\
    =&\mathcal{E}(u_m)+\mathcal{E}\Big(\sum_{j=1}^{q} \overline{U}_{p_j,\lambda_q}\Big)+o\Big(\frac{q}{\lambda_q^2}\Big)+O\Big(\frac{q}{\lambda_q^{\frac{n-2}{2}}}\Big)\nonumber\\
    =&\mathcal{E}(u_m)+q\Big(B_1+\frac{B_2}{\lambda_q^2}-\sum_{j=2}^{q}\frac{B_3}{\lambda_q^{n-2}|p_j-p_1|^{n-2}}\Big)+o\Big(\frac{q}{\lambda_q^2}\Big),
\end{align}
where $B_1, B_2, B_3$ are some positive constants.
And
\begin{equation}
    \frac{\partial F(\bar t_q,\bar z'_q,\lambda_q)}{\partial \lambda_q}= -q\Big(\frac{C_1}{\lambda_q^3}-\frac{C_2 q^{n-2}}{\lambda_q^{n-1}}+o\Big(\frac{1}{\lambda_q^{3}}\Big)\Big),
\end{equation}
where $C_1, C_2$ are the positive constants in Lemma \ref{lem4.7}.

In order to find a critical point for $F(\bar t_q,\bar z'_q,\lambda_q)$, we only need to make $c_l=0, l=3,\cdots,n-k$. Combining Proposition \ref{pro4.5} and Lemmas \ref{lem4.6}--\ref{lem4.7}, we conclude that there exists a $\rho\in(3\bar \delta, 4\bar \delta)$ such that the problem is equivalent to find a solution $(\bar t_q, \bar z'_q)$ of the following equations:
\begin{equation}\label{reee1}
    \frac{\partial K(\bar t,\bar z')}{\partial \bar t}=o\Big(\frac{1}{\lambda_q^{1/2}}\Big),
\end{equation}
\begin{equation}\label{reee2}
    \frac{\partial K(\bar t,\bar z')}{\partial \bar z_j}=o\Big(\frac{1}{\lambda_q^{1/2}}\Big),\quad j=5\cdots,n-k.
\end{equation}
\begin{equation}\label{reee3}
    C_1-C_2\frac{q^{n-2}}{\lambda_q^{n-4}}=o(1).
\end{equation}
Set $\lambda_q=\kappa q^{n-2/n-4}$,$\kappa\in[L_0, L_1]$, and
$$
\mathscr{G}(\kappa,\bar t_q, \bar z'_q):=\Big(\nabla_{\bar t_q, \bar z'_q}K(\bar t_q, \bar z'_q),C_1-\frac{C_2}{\kappa^{n-4}}\Big),
$$
then from ${\bf (K_1)}$ we have
\begin{equation*}
    \deg \Big(\mathscr{G}(\kappa,\bar t_q, \bar z'_q), [L_0, L_1]\times B_{\lambda_q^{\bar\theta-1}}(t_0, z'_0)\Big)= \deg \Big(\nabla_{\bar t_q, \bar z'_q}K(\bar t_q, \bar z'_q), B_{\lambda_q^{\bar\theta-1}}(t_0, z'_0)\Big)\neq 0.
\end{equation*}
Hence, \eqref{reee1}--\eqref{reee3} have a solution $(\bar t_q, \bar z'_q)$ satisfying
$|(\bar t_q, \bar z'_q)-(t_0,z'_0)|=o\Big(\frac{1}{\lambda_q^{1-\bar {\theta}}}\Big)$, and $\lambda_q\in[L_0 q^{\frac{n-2}{n-4}}, L_1q^{\frac{n-2}{n-4}}]$. Thus we have proved Theorem \ref{newsolution}.
\end{proof}

\medskip

\begin{remark}
    Our method of constructing new kinds of new cylindrial solutions can be applied to other kinds of critical Grushin problem. For example, the following equation with competing potentials:
    \begin{equation}\label{equationcompeting}
-\Delta u(x)+V(x)u(x)=\displaystyle K(x)\frac{u(x)^{2^*-1}}{|y|},\quad u>0 \quad \text{ in } \,  \R^n.
\end{equation}
Combining with the existence result in \cite{LiuNiuConstruction2022}, we can extend our existence result of new bubble solutions in Theorem \ref{newsolution} to \eqref{equationcompeting}. Since the idea of proof is very similar, in the following we give the statements of main results  for \eqref{equationcompeting} and leave the detailed proof for interested readers.

We assume  $V(x)$ satisfies:\\
${\bf (KV')}$:  $V(x)=V(|z^1|,z^2)\ge0$ and are bounded functions for $x=(y,z^1,z^2)\in \R^{k}\times\R^{2}\times\R^{n-k-2}$.  $V(r,z^2)\in C^1(B_{\rho_0}(r_0,z_0^2))$, $K(r,z^2)\in C^3(B_{\rho_0}(r_0,z_0^2))$
for $\rho_0>0$ is a fixed small constant, and
$$
V(r_0,z_0^2)\int_{\R^n}U_{0,1}^2\, \text{d}x-\frac{\Delta K(r_0,z_0^2)}{2^*(n-k)}\int_{\R^n}\frac{|z|^2}{|y|}U_{0,1}^{2^*}\, \text{d}x>0.
$$
We have the non-degeneracy result about the bubble solution, which we denote as $\widetilde {u}_m$, in \cite{LiuNiuConstruction2022}.

\begin{theorem}\label{nondeg'}
Suppose that $n\ge 8$, $\frac{n+1}{2}\le k<n-3$, $K(x)$ and $V(x)$ satisfies ${\bf (K_1')},{\bf (K_3')}$ and ${\bf (KV')}$, then there exists a large $\widetilde {m}_0$, such that for any integer $\widetilde {m}>\widetilde {m}_0$,
if $\varsigma\in {H_s}$ is a solution of the following linear equation:
$$
    \widetilde{L}_m \varsigma:=-\Delta\varsigma + V(x)\varsigma - (2^*-1)K(x)\frac{\widetilde {u}_m^{2^*-2}}{|y|}\varsigma=0\quad \text{ in } \, \R^n,
$$
then $\varsigma=0$.
\end{theorem}

Let  $z^*=(z_1,z_2,z_3,z_4)$ radially, we still denote the bubble solution centered at $\tilde \zeta_i$ as $\widetilde{u}_m$, and assume that:\\
${\bf (KV)}$:  $V(x)=V(|z^*|,z')\ge0$ and are bounded functions for $x=(y,z^*,z')\in \R^{k}\times\R^{4}\times\R^{n-k-4}$.  $V(t,z')\in C^1(B_{\rho_0}(t_0,z_0'))$, $K(t,z')\in C^3(B_{\rho_0}(t_0,z_0'))$
for $\rho_0>0$ is a fixed small constant, and
$$
V(t_0,z_0')\int_{\R^n}U_{0,1}^2\, \text{d}x-\frac{\Delta K(t_0,z_0')}{2^*(n-k)}\int_{\R^n}\frac{|z|^2}{|y|}U_{0,1}^{2^*}\, \text{d}x>0.
$$

As an application of the nondegeneracy result obtained in Theorem \eqref{nondeg'}, we have the following:
\begin{theorem}\label{newsolution'}
     Suppose that $n\ge 8$, $\frac{n+1}{2}\le k<n-3$, $K(x)$ and $V(x)$ satisfies ${\bf (K_1)},{\bf (K_3)}$ and ${\bf (KV)}$,  then there exists an integer $\widetilde {q_0}>0$, such that for any integer $\widetilde q>\widetilde {q_0}$, problem \eqref{equationcompeting} has a solution $\widetilde {v}_q$ of the form
\begin{equation*}\label{eq:newsolution'}
\widetilde {v}_q=\widetilde {u}_m+\sum_{j=1}^{\widetilde q} \eta U_{\widetilde {p}_j,\lambda_{\widetilde q}}+\widetilde {\psi}_q,
\end{equation*}
where $\widetilde {\psi}_q \in  {X_s}$, $(\bar t_{\widetilde q}, \bar z'_{\widetilde q}) \to (t_0, z'_0)$, $\lambda_{\widetilde q} \in [{\widetilde L_0} {\widetilde q}^{\frac{n-2}{n-4} }$, ${\widetilde L_1} {\widetilde q}^{\frac{n-2}{n-4}}],$ and $\|\widetilde {\psi}_q\|_{L^{\infty}(\R^n)}=o(\lambda_{\widetilde q}^{\frac{n-2}{2}})$.

Moreover, we deduce that the following Grushin problem with competing potentials for $K=K(|y|, z)={R(\sqrt{|y|}, z)}/{4},$ that is,
\begin{equation*}\label{eq:Grushin'}
    -\Delta_y u-4|y|^2 \Delta_z u+ 4|y|^2 V(y, z) u(y, z) =R(y, z) u(y, z)^{\frac{m_1+2 m_2+2}{m_1+2 m_2-2}},\quad(y, z) \in \mathbb{R}^{m_1} \times \mathbb{R}^{m_2},
\end{equation*}
has infinitely many cylindrically symmetric multi-bubbling solutions.
\end{theorem}
\end{remark}

\begin{remark}
    From the above theorems, we can conclude that $K(x)$ is the leader when competing with $V(x)$,  the bubble solutions only concentrate at the stable critical point $(t_0,z_0')$ of $K(x)$ and $V(x)$ has no affection on the non-degenerate condition ${\bf (K_3)}.$ The main reason of this phenomenon is because the related term $V(x)u(x)$ in \eqref{equationcompeting} usually decays faster than $K(x)\frac{u(x)^{2^*-1}}{|y|}$.
\end{remark}

\section*{Data availability}
No data was used for the research described in the article.

\medskip

\appendix
\section{Local Pohozaev identities}\label{appendixA}

This section is devoted to state the local Pohozaev identities for critical Hardy-Sobolev-type operator, which can be found in \cite{GuoMussoPengYanNon2020}. Let
$$
-\Delta u(x) =K(x)\frac{u^{2^*-1}(x)}{|y|},\quad u>0, \quad x=(y,z)\, \text{ in }\, \mathbb{R}^{k}\times\mathbb{R}^{n-k},
$$
and
$$
-\Delta \xi(x) =(2^*-1)K(x)\frac{u^{2^*-2}(x)}{|y|}\xi,\quad u>0, \quad x=(y,z)\, \text{ in }\, \mathbb{R}^{k}\times\mathbb{R}^{n-k}.
$$
Assume that $\Omega $ is a smooth bounded domain in $\mathbb{R}^{n}$. Then we have the following Lemma.
\begin{lemma}(Lemma 2.1, \cite{GuoMussoPengYanNon2020})\label{pohozeav}
It holds that
\begin{equation}\label{pohozeav1}
\begin{aligned}
      - \int_{\Omega}\frac{\partial K(r,z^2)}{\partial z_j}\frac{u^{2^*-1}\xi}{|y|}
      =-\int_{\partial\Omega}\Big(\frac{\partial u}{\partial \nu}\frac{\partial \xi}{\partial z_j}+\frac{\partial \xi}{\partial \nu}\frac{\partial u}{\partial z_j}\Big)+\int_{\partial\Omega}\nabla u\nabla \xi \nu_{k+j}-\int_{\partial\Omega}K(r,z^2)\frac{u^{2^*-1}\xi}{|y|}\nu_{k+j},
\end{aligned}
\end{equation}
and
\begin{equation}\label{pohozeav2}
\begin{aligned}
     &\int_{\Omega}\frac{u^{2^*-1} \xi }{|y|}\langle\nabla K(r,z^2),x-x_0\rangle\\=&\int_{\partial\Omega} \frac{K(r,z^2)}{|y|}u^{2^*-1}\xi \langle\nu,x-x_0\rangle
+\int_{\partial\Omega} \Big(\frac{\partial u}{\partial \nu}\langle\nabla \xi,x-x_0\rangle+\frac{\partial \xi}{\partial \nu}\langle\nabla u,x-x_0\rangle\Big)\\&-\int_{\partial\Omega} \nabla u \cdot\nabla\xi\langle\nu,x-x_0\rangle+\frac{n-2}{2}\int_{\partial\Omega}\Big(u\frac{\partial \xi}{\partial \nu}+\xi \frac{\partial u}{\partial \nu}\Big),
\end{aligned}
\end{equation}
where $j=1,\cdots,n-k$ and $\nu$ is the outer normal vector of $\Omega$.
\end{lemma}

\section{The Green's function}\label{appendixB}
In this part, we will establish the estimate of modified Green function, so that we obtain the  properties of the Green function of $L_m$, which is necessary for the construction of new cylindrical solutions. First, we need to define some corresponding operators.

Let $R_j$ as
$$
R_j x=\Big(y,\sqrt{z_1^2+z_2^2} \cos \big(\theta+\frac{2j\pi}{m}\big),\sqrt{z_1^2+z_2^2} \sin \big(\theta+\frac{2j\pi}{m}\big),z^2\Big),\quad j=1,\cdots,m,
$$
and let $T_i$  as
$$
T_i x =(y, z_1,\cdots,z_{i-1},(-1)\delta_{i2}z_i,z_{i+1},\cdots,z_n),\quad  i=1,\cdots,n-k,
$$
where $x=(y,z^1,z^2)\in \R^{k}\times\R^2\times\R^{n-k}.$
For any function $f$ defined in $\R^n$, define
$$
\bar f (y)= \frac{1}{m} \sum_{j=1}^m f(R_j y),
$$
and
$$
f^{*} (y)= \frac{1}{n-1}\sum_{i=2}^{n-k} \frac{1}{2}(\bar f(y)+\bar f(T_i y)).
$$
It is easy to check that $f^{*}\in H_s$.

To discuss the Green's function of $L_m$, regardless of $\delta_x$ not belonging to $H_s$, we consider
\begin{equation}\label{C1}
L_m u = \delta_x^{*}  \quad  \text{ in } \mathbb{R}^{n},\quad u \in H_s \cap \overline{D^{1,2}(\R^n) \cap H^1(\R^n)},
\end{equation}
where
$$
\delta_x^{*}=\frac{1}{n-k-1} \sum_{i=2}^{n-k}\frac{1}{2} \Big(\frac{1}{m} \sum_{j=1}^m \delta_{R_j x}+\frac{1}{m} \sum_{j=1}^m \delta_{T_i R_j x}\Big).
$$
We denote the solution of \eqref{C1} as $G_m(\tilde x,\bar x)$, which is called as the Green function of $L_m.$ We have

\begin{proposition}\label{proc.1}
The solution $G_m(\tilde x,\bar x)$ satisfies
\begin{equation}\label{gmc}
|G_m(\tilde x,\bar x)|\le \frac{C}{n-k-1} \sum_{i=2}^{n-k}\frac{1}{2} \Big(\frac{1}{m} \sum_{j=1}^m \frac{1}{|\tilde x-R_j \bar x|}+\frac{1}{m} \sum_{j=1}^m \frac{1}{|\tilde x-T_iR_j \bar x|}\Big)
\end{equation}
for all $\bar x\in B_R(0)$, where $R>0$ is any fixed large constant.
\end{proposition}

\begin{proof}
Let $v_1= G(\tilde x, x)$ be the Green's function of $-\Delta$ in $\R^n$. Let $v_2$ be the positive solution of
$$
\begin{cases}
\displaystyle -\Delta v= (2^*-1)K(r,z^2) \frac{u_m^{2^*-2}}{|y|} v_1 \quad &\text{ in }\, B_{2R}(0),\\
v=0\quad &\text{ on }\, \partial B_{2R}(0).
\end{cases}
$$
Then
$$
0\le v_2(\tilde x)\le \displaystyle(2^*-1)\int_{\R^n}G(\tilde x,  x) K(r,z^2) \frac{u_m^{2^*-2}}{|y|} v_1 \le C \frac{1}{|\tilde x-\bar x|^{n-3}}.
$$
We can continue this process to find $v_i$, which is the positive solution of
$$
\begin{cases}
-\Delta v = (2^*-1)K(r,z^2) \displaystyle\frac{u_m^{2^*-2}}{|y|} v_{i-1} \quad &\text{ in }\, B_{2R}(0),\\
v=0\quad &\text{ on }\, \partial B_{2R}(0).
\end{cases}
$$
And satisfies
$$
0\le v_i(\tilde x) \le C \frac{1}{|\tilde x-\bar x|^{n-1-i}}.
$$

Let $i$ be large enough so that $v_i\in L^\infty(B_{2R}(0))$.  Define
$$
v=\sum_{l=1}^i v_l\quad \text{ and }\quad w= G(\tilde x,\bar x)-  \iota v^*,
$$
where $\iota(x)\equiv\iota(z^1, z^2)\in C_0^{\infty}(B_{2R}(0))$, $\iota=1$ in $B_{\frac{3}{2}R}(0)$, and $0\le\iota\le1$.
Then we have
\begin{equation}\label{ccc}
\begin{cases}
L_k w = f \quad &\text{ in }\, B_{2R}(0),\\
w = 0\quad &\text{ on }\, \partial B_{2R}(0),
\end{cases}
\end{equation}
where $f\in L^\infty\cap {H_s}$.   By Theorem \ref{nondeg},  \eqref{ccc} has a solution $w\in {H}_s$.

By standard elliptic estimate, we have $w(\tilde x)$ is bounded, and
$$
|w(\tilde x)|\leq C\int_{\R^n}\frac{1}{|\tilde x-x|^{n-2}}\Big(\frac{|u_m(x)|^{2^*-2}}{|y|}|w(x)|+|g|\Big)\le\frac{C}{|\tilde x|}.
$$
Then we can continue this process and finally prove \eqref{gmc}.
\end{proof}

\medskip
\section{Basic estimates and lemmas}\label{appendixC}
This section is devoted to state some useful and well-known estimates and lemmas.

\begin{lemma}\label{B.1}
Assume that $\alpha> 0,$ we have the following estimates for $m \to +\infty$, $j=2,\cdots,m$:
\begin{equation}\label{b1}
\begin{aligned}
\sum_{j=2}^m  \frac{1}{|\zeta_{1}-\zeta_{j}|^{\alpha}}
& =
\begin{cases}
\displaystyle O\Big(\frac{m^\alpha}{\bar r^{\alpha}}\Big) & \text { if }\, \alpha>1 ,\\
\displaystyle O\Big(\frac{m^\alpha \ln m}{\bar r^{\alpha}}\Big) & \text { if }\, \alpha=1,\\
\displaystyle O\Big(\frac{m }{\bar r^{\alpha}}\Big) & \text { if }\, \alpha<1.
\end{cases}
\end{aligned}
\end{equation}
\end{lemma}

\begin{proof}
The proof of Lemma \ref{B.1} is similar to that of Lemma A.3 in \cite{GaoGuoNewarXiv}, here we omit it.
\end{proof}

Define
$$
g_{ij}(y)=\frac{1}{(1+|y|+|z-\zeta_{i}|)^{\gamma_{1}}}\frac{1}{(1+|y|+|z-\zeta_{j}|)^{\gamma_{2}}},\quad i \neq j,
$$
where  $\gamma_{1}\geq 1 $ and  $\gamma_{2}\geq 1 $ are two constants.

\begin{lemma}\label{B.5}(Lemma A.1, \cite{LiuWangCylindrical2023})
For any constants  $0<\upsilon\leq \min\{\gamma_{1},\gamma_{2}\} $, there is a constant  $C>0 $, such that
$$
g_{ij}(y)\leq \frac{C}{|\zeta_{i}-\zeta_{j}|^{\upsilon}}\Big(\frac{1}{(1+|y|+|z-\zeta_{i}|)^{\gamma_{1}+\gamma_{2}-\upsilon}}+\frac{1}{(1+|y|+|z-\zeta_{j}|)^{\gamma_{1}+\gamma_{2}-\upsilon}}\Big).
$$
\end{lemma}

\begin{lemma}\label{B.6}(Lemma A.2, \cite{LiuWangCylindrical2023})
Assume that $n\ge 5$, $\frac{n+1}{2}\le k <n-1$. Then for any constant  $0<\beta <n-2$, there is a constant  $C>0$, such that for all $x=(y,z)\in \R^{k}\times\R^{n-k}$,
$$
\int_{\R^n} \frac{1}{|\tilde x-x|^{n-2}}\frac{1}{|\tilde y|(1+|\tilde y|+|\tilde z-\zeta_i|)^{1+\beta}}\, {\mathrm d}\tilde x \leq
\frac{C}{(1+|y|+|z-\zeta_i|)^{\beta}}.
$$
\end{lemma}

\begin{lemma}\label{B.7}(Lemma A.3, \cite{LiuWangCylindrical2023})
Assume that $n\ge 5$, $\frac{n+1}{2}\le k <n-1$. Then there is a constant $C>0$ and a small
$\theta>0$, such that for all $x=(y,z)\in \R^{k}\times\R^{n-k}$,
$$
\int_{\R^n}\frac1{|\tilde x-x|^{n-2}} \frac{\overline{W}_{r,h,\mu}^{2^*-2}(\tilde x)}{|\tilde y|}\sum_{j=1}^k\frac1{(1+|\tilde y|+|\tilde z-\zeta_i|)^{\frac{n-2}{2}+\tau}}\, {\mathrm d}\tilde x \le C\sum_{j=1}^k\frac1{(1+|y|+|z-\zeta_i|)^{\frac{n-2}{2}+\tau+\theta}},
$$
where $\tau=\frac{n-4}{n-2}$.
\end{lemma}

\end{document}